\documentclass{article}

\usepackage[utf8]{inputenc}
\usepackage[T1]{fontenc}
\usepackage{refcheck}
\usepackage{fullpage}
\usepackage{xargs}
\usepackage{amsmath,amssymb,fancybox,pgf,array,amsthm}

\renewcommand{\b}[1]{{\bf #1}}

\newcommand{\fz}[3]{#1:\, #2 \rightarrow #3}

\renewcommand{\r}[1]{(\ref{#1})}
\newcommand{\bi}{\begin{itemize}}
\newcommand{\ei}{\end{itemize}}
\newcommand{\be}{\begin{enumerate}}
\newcommand{\ee}{\end{enumerate}}

\newcommand{\bd}{\begin{description}}
\newcommand{\ed}{\end{description}}
\renewcommand{\i}{\item}

\newcommand{\bqn}{\begin{eqnarray}}
\newcommand{\eqn}{\end{eqnarray}}
\newcommand{\eqnn}{\nonumber\end{eqnarray}}
\newcommand{\eqnl}[1]{\label{#1}\end{eqnarray}}

\newcommand{\nn}{\nonumber\\}

\newcommand{\ba}[1]{\begin{array}{#1}}
\newcommand{\ea}{\end{array}}

\newcommand{\R}{\mathbb{R}}

\newcommand{\N}{\mathbb{N}}

\newcommand{\bproof}{\begin{proof}}
\newcommand{\eproof}{\end{proof}}
\newtheorem{Theorem}{\bf Theorem}
\newtheorem{lemma}[Theorem]{\bf Lemma}
\newtheorem{corollary}[Theorem]{\bf Corollary}
\newtheorem{definition}[Theorem]{\bf Definition}
\newtheorem{proposition}[Theorem]{\bf Proposition}
\newtheorem{remark}[Theorem]{\bf Remark}

\newcommand{\bt}{\begin{Theorem}}
\newcommand{\et}{\end{Theorem}}
\newcommand{\bl}{\begin{lemma}}
\newcommand{\el}{\end{lemma}}
\newcommand{\bp}{\begin{proposition}}
\newcommand{\ep}{\end{proposition}}
\newcommand{\bc}{\begin{corollary}}
\newcommand{\ec}{\end{corollary}}
\newcommand{\bdeff}{\begin{definition}}
\newcommand{\edeff}{\end{definition}}
\newcommand{\brem}{\begin{remark}\rm}
\newcommand{\erem}{\end{remark}}
\newcommand{\auth}[1]{{\sc #1}}
\newcommand{\tit}[1]{{\rm #1}}
\newcommand{\titl}[1]{{\it #1}}
\newcommand{\jou}[1]{{\it #1}}

\newcommand{\pp}[1]{pp.~#1}

\newcommand{\lam}{\lambda}
\newcommand{\al}{\alpha}
\newcommand{\eps}{\varepsilon}
\newcommand{\ga}{\gamma}

\newcommand{\Id}{\mathrm{Id}}

\newcommand{\Pt}[1]{\left( #1 \right)}
\newcommand{\Pg}[1]{\left\{ #1 \right\}}
\newcommand{\Pq}[1]{\left[ #1 \right] }
\newcommand{\Pa}[1]{\langle #1 \rangle}
\newcommand{\Pabs}[1]{{\Big \vert}  #1 {\Big \vert}}

\newcommand{\Mu}{\mathcal{M}}

\newcommand{\weak}{\rightharpoonup}
\newcommand{\supp}{\mathrm{supp}}
\newcommand{\sotto}{\leq}
\newcommand{\Der}[2]{{D_{#2}#1}}

\begin{document}

\newcommand{\gwbase}{W^{a,b}_p}
\newcommand{\gw}[2][p]{W^{a,b}_{#1}\Pt{#2}}
\newcommand{\wwbase}{\mathcal{W}}
\newcommand{\ww}[1]{\mathcal{W}\Pt{#1}}

\newcommand{\tw}[2][p]{T^{a,b}_{#1}\Pt{#2}}
\newcommand{\twbase}{T^{a,b}_p}

\newcommand{\LL}{\mathbb{L}}
\newcommand{\KK}{\mathbb{K}}
\newcommand{\FF}{\mathbb{F}}
\newcommand{\DD}{\mathcal{D}}

\newcommand{\Muu}{{\Mu_{0}^{ac}}}
\renewcommand{\P}{\mathcal{P}}
\newcommand{\Pac}{{\P_{0}^{ac}}}
\newcommand{\inter}{\mathrm{int}}

\newcommand{\gwbbase}{W^{a,b}_2}
\newcommand{\gwbb}[1]{\gwbbase\Pt{#1}}
\newcommand{\fa}[1]{\mathcal{A}\Pq{#1}}
\newcommand{\fabase}{\mathcal{A}}
\newcommand{\fb}[1]{\mathcal{B}^{a,b}\Pq{#1}}
\newcommand{\fbbase}{\mathcal{B}^{a,b}}
\newcommand{\bart}{{\bar t}}
\newcommand{\barx}{{\bar x}}
\newcommand{\bmu}[1][k]{\mu^{\Pq{#1}}}
\newcommand{\tmu}[1][k]{\tilde{\mu}^{\Pq{#1}}}
\newcommand{\ff}[2][k]{\Phi^{#1}_{\Pq{#2}}}
\newcommand{\hhm}[2][k]{\underline{H}^{#1}_{\Pq{#2}}}
\newcommand{\hhp}[2][k]{\overline{H}^{#1}_{\Pq{#2}}}
\newcommand{\bbmu}{\tilde{\mu}^{[k]}}
\newcommand{\bbv}{\tilde{v}^k}
\newcommand{\bbh}{\tilde{h}^k}
\newcommand{\dt}{{\Delta t}}

\newcommand{\pho}{\rho}
\newcommand{\bbar}[1]{{\bar{\bar{#1}}}}

\newcommand{\wubase}{W_1^{1,1}}
\newcommand{\wu}[1]{W_1^{1,1}\Pt{#1}}
\newcommand{\nuovo}[1]{{\bf #1}}

\title{On properties of the Generalized Wasserstein distance}


\author{Benedetto Piccoli\thanks{Department of Mathematical Sciences, Rutgers University - Camden, Camden, NJ. {\tt piccoli@camden.rutgers.edu}}, Francesco Rossi\thanks{Aix Marseille Universit\'e, CNRS, ENSAM, Universit\'e de Toulon, LSIS UMR 7296, 13397, Marseille, France {\tt francesco.rossi@lsis.org}}}

\maketitle

\begin{abstract}

The Wasserstein distances $W_p$ ($p\geq 1$), defined in terms of solution to the Monge-Kantorovich problem, are known to be a useful tool
to investigate transport equations.
In particular, the Benamou-Brenier formula characterizes the
square of the Wasserstein distance $W_2$ as the infimum
of the kinetic energy, or action functional, of all vector fields 
moving one measure to the other.\\
Another important property of the Wasserstein distances is the
Kantorovich-Rubinstein duality stating the equality between the
distance $W_1$ and the supremum of the integrals of 
Lipschitz continuous functions with Lipschitz constant bounded by one.

An intrinsic limitation of Wasserstein distances is the fact
that they are defined only between measures having the same mass.
To overcome such limitation, we recently introduced 
the generalized Wasserstein distances $W_p^{a,b}$,
defined in terms of both the classical Wasserstein distance $W_p$
and the total variation (or $L^1$) distance, see \cite{gw}.
Here $p$ plays the same role as for the classic Wasserstein distance,
while $a$ and $b$ are weights for the transport and the total variation term.\\
In this paper we prove two important properties of the generalized
Wasserstein distances:\\
1) a generalized Benamou-Brenier formula
providing the equality between $W_2^{a,b}$ and the supremum of an action
functional, which includes a transport term (kinetic energy) and
a source term.\\
2) a duality \emph{\`{a} la} Kantorovich-Rubinstein
establishing the equality
between $W_1^{1,1}$ and the flat metric.\\
\end{abstract}



\section{Introduction}

The problem of optimal transportation, also called  Monge-Kantorovich problem, has been intensively studied in mathematical community. Related to this problem, Wasserstein distances in the space of probability measures have revealed to be powerful tools, in particular for dealing with dynamics of measures (like the transport PDE, see e.g. \cite{ambrosio,gradient}). For a complete introduction to Wasserstein distances, see \cite{old-new,villani}.

The main limit of this approach, at least for its application to dynamics of measures, is that the Wasserstein distances $W_p(\mu,\nu)$ ($p\geq 1$)
are defined only if the two measures $\mu,\nu$ have the same mass. For this reason, in \cite{gw} we introduced the generalized Wasserstein distances $\gw[p]{\mu,\nu}$, combining the standard Wasserstein and 
total variation distances. In rough words, for $\gw[p]{\mu,\nu}$ an infinitesimal mass $\delta\mu$ of $\mu$ can either be removed at cost $a |\delta\mu|$, or moved from $\mu$ to $\nu$ at cost $bW_p(\delta\mu,\delta\nu)$. More formally, the definition of the generalized Wasserstein distance that we use in this article\footnote{Observe that the definition in \cite{gw} was $\gw{\mu,\nu}=\inf_{\tilde\mu,\tilde\nu\in\Mu} a|\mu-\tilde\mu|+a|\nu-\tilde\nu|+bW_p(\tilde\mu,\tilde\nu)$. Clearly, the two definitions are extremely similar, and satisfy similar properties: one can indeed observe that , given the vector 
$(a|\mu-\tilde\mu|+a|\nu-\tilde\nu|,bW_p(\tilde\mu,\tilde\nu))\in\R^2$, the definition in \cite{gw} is the 1-norm of such vector, while the definition given in the present article is its $p$-norm.} is
\bqn
\gw{\mu,\nu}:=\Pt{\tw{\mu,\nu}}^{1/p},
\eqnn
with
\bqn
\tw{\mu,\nu}=\inf_{\tilde\mu,\tilde\nu\in\Mu,\,|\tilde\mu|=|\tilde\nu|}a^p\Pt{|\mu-\tilde\mu|+|\nu-\tilde\nu|}^p+b^pW_p^p(\tilde\mu,\tilde\nu),
\eqnn
where $\Mu$ denotes the space of Borel regular measures on $\R^d$ with finite mass.

Recall that the ``flat metric'' or ``bounded Lipschitz distance'' (see e.g. \cite[\S 11]{dudley}), is defined as follows
$$
d(\mu,\nu):=\sup\Pg{\int_{\R^d} f\,d(\mu-\nu)\,|\, \|f\|_{C^0}\leq 1, \|f\|_{Lip}\leq 1}.
$$
We first show that the generalized Wasserstein distance $W_1^{1,1}$ coincides with the flat metric. This provides the following duality formula: $$d(\mu,\nu)=W_1^{1,1}(\mu,\nu)=
\inf_{\tilde\mu,\tilde\nu\in\Mu,\,|\tilde\mu|=|\tilde\nu|} |\mu-\tilde\mu|+|\nu-\tilde\nu|+W_1(\tilde\mu,\tilde\nu).$$
This result can be seen as a generalization of the  Kantorovich-Rubinstein theorem, which provides the duality:
$$
W_1(\mu,\nu)=\sup\Pg{\int_{\R^d} f\,d(\mu-\nu)\,|\,  \|f\|_{Lip}\leq 1}.
$$\\

One interesting field of application of the generalized Wasserstein distances is the study of transport equations with sources, i.e. dynamics of measures given by:
\bqn
\partial_t\mu_t+\nabla\cdot(v_t \mu_t)=h_t,
\eqnl{pde}
where $v_t$ is a time-dependent vector field and $h_t$ a time-dependent source term.
Several authors have studied \r{pde} without source term, i.e. $h\equiv 0$, showing that it is very convenient to use the standard Wasserstein distance in this framework.  In particular, Benamou and Brenier showed in \cite{ben-bren} that there is a natural equivalence between the minimization of the action functional $\fa{\mu,v}:=\int_0^1 dt\Pt{\int_{\R^d}d\mu_t |v_t|^2}$ and the computation of the Wasserstein distance $W_2$. Their fundamental result is recalled in Theorem \ref{t-bb}.
However, the standard Wasserstein distances do not encompass the case of a non vanishing source $h$. Indeed, in this case the mass of the measure $\mu_t$ varies in time, hence $W_p(\mu_t,\mu_s)$ may not be defined for $t\neq s$.\\
Our second goal is to generalize the Benamou-Brenier formula to this setting. On one side, we use the generalized Wasserstein distances, so allowing mixing creation/removal of mass and transport of mass. On the other side, we define a generalization of the functional $\fabase$, taking into account both the transport and the creation/removal of mass in \r{pde}. More precisely, we define
$$\fb{\mu,v,h}:=a^2\Pt{\int_0^1 dt \Pt{\int_{\R^d}d |h_t|}}^2 +b^2\int_0^1 dt \Pt{ \int_{\R^d} d\mu_t\,|v_t|^2}.$$
Given the generalizations both for the distance and the functional, we will then prove the generalized Benamou-Brenier formula under the regularity hypotheses recalled in Definition \ref{d-a}:
\bqn
\tw[2]{\mu_0,\mu_1}=\inf\Pg{\fb{\mu,v,h}\ \Big|\ 
\ba{c}
\mu\mbox{~ is a solution of \r{pde} with vector field $v$, source $h$}\\
\mbox{and }\mu_{|_{t=0}}=\mu_0,\mu_{|_{t=1}}=\mu_1
\ea}.
\eqnl{intro-bb}\\

The structure of the paper is the following. In Section \ref{s-gw} we define the generalized Wasserstein distance and recall some useful properties, in particular estimates of the generalized Wasserstein distance under flow action. In Section \ref{s-flat} we prove that $\wubase$ coincides with the flat metric. Finally, in Section \ref{s-bb} we recall the standard Benamou-Brenier formula and prove the generalized Benamou-Brenier formula \r{intro-bb}.
\section{Generalized Wasserstein distance}
\label{s-gw}


\subsection{Notation and standard Wasserstein distance}
\label{s-monge}


We use $\Mu$ to denote the space of positive Borel regular measures with finite mass\footnote{The requirement of having finite mass is a simple choice to have finite distances $\gw{\mu,nu}$.} on $\R^d$ and $\Muu$ to denote the subspace of $\Mu$ of measures with compact support that are absolutely continuous with respect to the Lebesgue measure. \nuovo{When not specified, the domain of integration is the whole space $\R^d$, or $\R^d\times \R^d$ in the case of integrals with two variables.}

Given $\mu,\mu_1$ Radon measures (i.e. positive Borel measures with locally finite mass), we write $\mu_1\ll \mu$
if $\mu_1$ is absolutely continuous with respect to $\mu$, while we write $\mu_1\sotto\mu$
if $\mu_1(A)\leq\mu(A)$ for every Borel set $A$. We denote with $|\mu|:=\mu(\R^d)$ the norm 
of $\mu$ (also called its mass). More generally, if $\mu=\mu^+-\mu^-$ is a signed Borel measure, we define $|\mu|:=|\mu^+|+|\mu^-|$. 

By the Lebesgue's decomposition theorem, given two measures $\mu,\nu$, one can always write in a unique way $\mu=\mu_{ac}+\mu_s$ such that $\mu_{ac}\ll \nu$ and $\mu_{s}\perp\nu$, i.e. there exists $B$ such that $\mu_s(B)=0$ and $\nu(\R^n\setminus B)=0$. Moreover,
there exists a unique $f\in L^1(d\nu)$ such that $d\mu_{ac}(x)=f(x)\,d\nu(x)$. Such $f$ is called the Radon-Nikodym derivative of $\mu$ with respect to $\nu$. We denote it with $\Der{\mu}{\nu}$. For more details, see e.g. \cite{ev-gar}.\\

Given a Borel map $\fz{\gamma}{\R^d}{\R^d}$, the push forward of a measure $\mu\in\Mu$
is defined by:
\bqn
\gamma\#\mu(A):=\mu(\gamma^{-1}(A)).
\eqnn
Note that the mass of $\mu$ is identical to the mass of $\gamma\#\mu$.
Therefore, given two measures $\mu,\nu$ with the same mass, one may look for $\gamma$
such that  $\nu=\gamma\#\mu$ and it minimizes the cost
$$I\Pq{\gamma}:=|\mu|^{-1}\,\int |x-\ga(x)|^p \,d\mu(x).$$
This means that each infinitesimal mass $\delta\mu$ is sent to $\delta \nu$ and that its infinitesimal cost is related to the $p$-th power of the distance between them. 
Such minimization problem is known as the Monge problem and was first stated by 1781 (see \cite{monge}).\\
If $\mu$ or $\nu$ has an atomic part then we may have no $\gamma$ such that
$\gamma\#\mu$. For instance, $\mu=2\delta_1$ and $\nu=\delta_{0}+\delta_2$, measures on the real line, have the same mass, but there exists no $\gamma$ with $\nu=\gamma\#\mu$, since $\gamma$ cannot separate masses. 
A simple condition, that ensures the existence of a minimizing $\gamma$, is that $\mu$ and $\nu$ are absolutely continuous with respect to the Lebesgue measure.

A generalization of the Monge problem is achieved as follows.
Given a probability measure $\pi$ on $\R^d\times\R^d$, one can interpret $\pi$ as a method to transfer a measure $\mu$ on $\R^d$ to another measure on $\R^d$ as follows: each infinitesimal mass on a location $x$ is sent to a location $y$ with a probability given by $\pi(x,y)$. Formally, $\mu$ is sent to $\nu$ if the following properties hold:
\bqn
|\mu|\,\int_{\R^d} d\pi(x,\cdot)=d\mu(x),\qquad \qquad |\nu|\,\int_{\R^d} d\pi(\cdot,y)=d\nu(y).
\eqnl{e-pi}
Such $\pi$ is called a transference plan from $\mu$ to $\nu$. We denote the set of such transference plans as $\Pi(\mu,\nu)$. Since one usually deals with probability measures $\mu,\nu$, the terms $|\mu|,|\nu|$ are usually neglected in the literature. A condition equivalent to \r{e-pi} is that, for all $f,g\in C^\infty_c(\R^d)$ it holds $|\mu|\,\int_{\R^d\times \R^d} (f(x)+g(y))\,d\pi(x,y) = \int_{\R^d} f(x)\,d\mu(x)+ \int_{\R^d} g(y)\,d\nu(y)$.
\brem \label{r-sotto} One can use a transference plan $\pi\in \Pi(\mu,\nu)$ also to define pairs $\mu'\sotto \mu,\nu'\sotto\nu$ so that $\mu'$ is transfered to $\nu'$. Indeed, given $\mu'\sotto \mu$, define the Radon-Nikodym derivative $f=\Der{\mu'}{\mu}$, that satisfies $f\leq 1$ and $\mu'(A)=\int_A f(x) d\mu(x)$ for all Borel sets. Define now $\pi',\nu'$ as follows:
\bqn
\pi'(A\times B)&:=&\frac{|\mu|}{|\mu'|} \int_{A\times B}f(x)d\pi(x,y)\mbox{~~~ for each Borel set~~}A\times B,\nn
\nu'(B)&:=&|\mu|\int_{\R^d\times\Pg{B}}f(x)d\pi(x,y)\mbox{~~~ for each Borel set~~}B.
\eqnn
It is easy to prove that $\pi'\in\Pi(\mu',\nu')$. Similarly, one can define $\pi''\in\Pi(\mu-\mu',\nu-\nu')$ by 
\bqn
\pi''(A\times B)&:=&\frac{|\mu|}{|\mu|-|\mu'|}\int_{A\times B}(1-f(x))d\pi(x,y)\mbox{~~~ for each Borel set~~}A\times B.
\eqnn

By semplicity, we will drop the passage from $\pi$ to $\pi'$ from now on. We will say that, given a transference plan $\pi\in\Pi(\mu,\nu)$ and $\mu'\sotto\mu$, then there exists a unique $\nu'$ such that $\pi\in\Pi(\mu',\nu')$.
\erem

One can define a cost for $\pi$ as follows $$J\Pq{\pi}:=\int_{\R^d\times\R^d} |x-y|^p \,d\pi(x,y)$$ and look for a minimizer of $J$ in $\Pi(\mu,\nu)$. Such problem is called the Monge-Kantorovich problem. It is important to observe that such problem is a generalization of the Monge problem. Indeed, given a $\gamma$ sending $\mu$ to $\nu$, one can define a transference plan $\pi=(\Id\times \gamma)\#\mu$, i.e. $d\pi(x,y)=\mu(\R^n)^{-1}\, d\mu(x)\delta_{y=\gamma(x)}$. It also holds $J\Pq{\Id\times\gamma}=I\Pq{\gamma}$. The main advantage of this approach is that a minimizer of $J$ in $\Pi(\mu,\nu)$ always exists.\\

A natural space on which $J$ is finite is the space of Borel measures with finite $p$-moment, that is
\bqn
\Mu^p:=\Pg{\mu\in\Mu\ |\ \int |x|^p\, d\mu(x)<\infty}. 
\eqnn
One can thus define on $\Mu^p$ the following operator between measures of the same mass\footnote{Remark that in \cite{gw} we hade the mass coefficient $|\mu|^{1/p}$. The choice here helps to have estimates not depending on $p$.}, called the \b{Wasserstein distance}:
\bqn
W_p(\mu,\nu)=|\mu|(\min_{\pi\in\Pi(\mu,\nu)} J\Pq{\pi})^{1/p}.
\eqnn
It is indeed a distance on the subspace of measures in $\Mu^p$ with a given mass, see \cite{villani}. It is easy to prove that $W_p(k\mu,k\nu)=k W_p(\mu,\nu)$ for $k\geq 0$, by observing that $\Pi(k\mu,k\nu)=\Pi(\mu,\nu)$ and that $J\Pq{\pi}$ does not depend on the mass.

Another remarkable property is the following principle for optimality.
\bp \label{p-splitWp}
Let $\pi\in\Pi(\mu,\nu)$ be a transference plan realizing $W_p(\mu,\nu)$. Let $\mu'\sotto\mu$ and $\nu'\sotto\nu$ such that $\pi\in\Pi(\mu',\nu')$. Then $\pi$ also realizes $W_p(\mu',\nu')$ and it holds
\bqn
\frac{W_p^p(\mu,\nu)}{|\mu|^{p-1}}=
\frac{W_p^p(\mu',\nu')}{|\mu'|^{p-1}}+
\frac{W_p^p(\mu-\mu',\nu-\nu')}{|\mu-\mu'|^{p-1}}.
\eqnl{e-split}
\ep
\bproof
First observe the precise meaning of the statement: define $\pi'$ the restriction of $\pi$ to $\mu',\nu'$ and with $\pi''$ the restriction of $\pi$ to $\mu-\mu',\mu-\nu'$, as explained in Remark \ref{r-sotto}. Then $\pi'$ is the transference plan realizing  $W_p(\mu',\nu')$. Also observe that $|\mu| J[\pi]=|\mu'| J[\pi']+ (|\mu|-|\mu'|) J[\pi'']$.

We first prove that $\pi'$ realizes $W_p(\mu',\nu')$, by contradiction. Assume that there exists $\tilde \pi'\in\Pi(\mu',\nu')$ such that $J\Pq{\tilde\pi'}<J\Pq{\pi'}$. Then define the transference plan $\tilde\pi\in\Pi(\mu,\nu)$ as follows:
\bqn
\tilde\pi(A\times B)&:=&\frac{|\mu'|}{|\mu|} \tilde\pi'(A\times B)+ \frac{|\mu|-|\mu'|}{|\mu|} \pi''(A\times B)\mbox{~~~ for each Borel set~~}A\times B.
\eqnn
A direct computation shows that
\bqn
|\mu| J[\tilde\pi]=|\mu'| J[\tilde\pi']+ (|\mu|-|\mu'|) J[\pi'']<|\mu'| J[\pi']+ (|\mu|-|\mu'|) J[\pi'']=|\mu| J[\pi].
\eqnn
Then $J[\tilde\pi]<J[\pi]$ and $\tilde\pi\in\Pi(\mu,\nu)$. This is in contradiction with the fact that $\pi$ realizes $W_p(\mu,\nu)$.

We have just proved that $\pi'$ realizes $W_p(\mu',\nu')$. By symmetry, we also have that $\pi''$ realizes $W_p(\mu-\mu',\nu-\nu')$. Then, the proof of \r{e-split} is a direct consequence of the fact that $|\mu| J[\pi]=|\mu'| J[\pi']+ (|\mu|-|\mu'|) J[\pi'']$.
\eproof

\subsection{Definition of the generalized Wasserstein distance}

In this section, we provide a definition of the generalized Wasserstein distance, which is a slight modification of that given in \cite{gw}, together with some useful properties.

\bdeff
Let $\mu,\nu\in\Mu$ be two measures. We define the functionals
\bqn
\tw{\mu,\nu}:=\inf_{\tilde\mu,\tilde\nu\in\Mu,\,|\tilde\mu|=|\tilde\nu|}a^p\Pt{|\mu-\tilde\mu|+|\nu-\tilde\nu|}^p+b^pW_p^p(\tilde\mu,\tilde\nu),
\eqnl{e-tw}
and
\bqn
\gw{\mu,\nu}:=\Pt{\tw{\mu,\nu}}^{1/p}.
\eqnl{e-gw}
\edeff

We now provide some properties of $\gwbase$ and $\twbase$. Proofs can be adapted from those given in \cite{gw}. 
\bp \label{p-base}
The following properties hold:\\
1. The infimum in \r{e-tw} coincides with
$$\inf_{\tilde\mu\sotto\mu,\tilde\nu\sotto\nu,\,|\tilde\mu|=|\tilde\nu|}a^p\Pt{|\mu-\tilde\mu|+|\nu-\tilde\nu|}^p+b^pW_p^p(\tilde\mu,\tilde\nu),$$
where we have added the constraint $\tilde\mu\sotto\mu,\tilde\nu\sotto\nu$.\\
2. The infimum in \r{e-tw} is attained by some $\tilde\mu,\tilde\nu$.\\
3. The functional $\gwbase$ is a distance on $\Mu$.\\
4. It holds $\gw{\mu,0}\leq a|\mu|$
\ep

\brem One could define another metric, similar to $\gwbase$, by replacing $\twbase$ with
$$\inf_{\tilde\mu,\tilde\nu\in\Mu,\,|\tilde\mu|=|\tilde\nu|}a^p\Pt{|\mu-\tilde\mu|^p+|\nu-\tilde\nu|^p}+b^pW_p^p(\tilde\mu,\tilde\nu),$$
i.e. by distributing the $p$-th power on the two $L_1$ terms. Proofs and properties are similar to the proofs given here.
Our choice here is related to the generalization of the Benamou-Brenier formula for $\gwbase$. We discuss this issue in Remark \ref{r-altrib} below.
\erem

We also have this useful estimate to bound integrals.
\bl
Let $\mu,\nu\in \Muu$, and $f\in \mathrm{Lip}(\R^d,\R)\cap L^\infty (\R^d,\R)$. Then
\bqn
\Pabs{\int f\, d\mu-\int f \, d\nu}\leq \sqrt{2}\max\Pg{\frac{\|f\|_\infty}{a},\frac{\|f\|_{Lip}}{b}}\gw{\mu,\nu}.
\eqnl{e-gwintegrale}
\el
\bproof
Let $\tilde\mu\sotto\mu,\tilde\nu\sotto\nu$ realizing $\gw{\mu,\nu}$. We have 
\bqn
\Pabs{\int f\, d\mu-\int f \, d\nu}&\leq& \Pabs{\int f\, d(\mu-\tilde\mu)}+\Pabs{\int f\, d(\tilde \mu-\tilde\nu)}+
\Pabs{\int f\, d(\tilde\nu-\nu)}\leq\nn
&\leq& \|f\|_\infty |\mu-\tilde\mu|+\|f\|_{Lip} W_1(\tilde\mu,\tilde\nu)+ \|f\|_\infty |\tilde\nu-\nu|,
\eqnl{e-ovvio}
where we have used that $|\mu|=\sup\Pg{\int f \,d\mu\ |\ \|f\|_\infty=1}$ and the Kantorovich-Rubinstein duality formula 
$W_1(\mu,\nu)=\sup\Pg{\int f \,d(\mu-\nu)\ |\ \|f\|_{Lip}=1}$. Recall that $W_1(\tilde\mu,\tilde\nu)\leq W_p(\tilde\mu,\tilde\nu)$ for $p\geq 1$, see e.g. \cite[Sec. 7.1.2]{villani}. Then \r{e-gwintegrale} is a direct consequence of \r{e-ovvio}, by using $(x+y)^2\leq 2(x^2+y^2)$.
\eproof
\subsection{Topology of the generalized Wasserstein distance}
In this section we recall some useful topological results related to the metric space $\Mu$ when endowed with the generalized Wasserstein distance.
We first define tightness in this context.
\bdeff \label{d-tight} A set of measures $M$ is tight if for each $\eps>0$ there exists a compact $K_\eps$ such that $\mu(\R^d\setminus K_\eps)<\eps$ for all $\mu\in M$.
\edeff

We now recall the following important result about convergence with respect to the generalized Wasserstein distance, see \cite[Theorem 13]{gw}.
\bt \label{t-convergence}
Let $\Pg{\mu_n}$ be a sequence of measures in $\R^d$, and $\mu_n,\mu\in\Mu$. Then
$$\gw{\mu_n,\mu} \to 0\mbox{~~~~~~is equivalent to~~~~~~}\mu_n\weak \mu \mbox{~~and~~}\Pg{\mu_n}\,\mbox{is tight}.$$
\et

We finally recall the result of completeness, see \cite[Proposition 15]{gw}.
\bp \label{p-complete}
The space $\Mu$ endowed with the distance $\gwbase$ is a complete metric space.
\ep

\subsection{Estimates of generalized Wasserstein distance under flow actions}
\label{s-flow}

In this section we give useful estimates both for the standard and generalized Wasserstein distances $W_p$ and $\gwbase$ under flow actions. Similar\footnote{Properties proven in \cite{nostro,gw} were not optimal, since we had a coefficient $e^{\frac{p+1}{p}Lt}$ instead of the coefficient $e^{Lt}$ in properties 1 and 3, and a coefficient $e^{Lt/p}$ instead of 1 in property 3.} properties were already proved for measures $\mu,\nu\in\Muu$ in \cite[Sec. 2.1]{nostro} and \cite[Sec.1.5]{gw}. Generalizations of these estimates to any measures in $\Mu$ are obvious, by using the Kantorovich formulation of the optimal transportation problem.

\bp \label{p-flow} Let $v_t,w_t$ be two time-varying vector fields, uniformly Lipschitz with respect to the space variable, and $\phi^t,\psi^t$ the flow generated by $v,w$ respectively. Let $L$ be the Lipschitz constant of $v$ and $w$, i.e. $|v_t(x)-v_t(y)|\leq L |x-y|$ for all $t$, and similarly for $w$. Let $\mu,\nu\in\Mu$. We have the following estimates for the standard Wasserstein distance
\bi
\i $W_p\Pt{\phi^t\#\mu,\phi^t\#\nu}\leq e^{Lt}W_p\Pt{\mu,\nu}$,
\i $W_p\Pt{\mu,\phi^t\#\mu}\leq t\|v\|_{C^0} |\mu|$,
\i $W_p\Pt{\phi^t\#\mu,\psi^t\#\nu}\leq e^{Lt}W_p\Pt{\mu,\nu}+\frac{e^{Lt}-1}{L}\,|\mu|\,\sup_{\tau\in[0,t]}\|v_t-w_t\|_{C^0} $.
\ei
We have the following estimates for the generalized Wasserstein distance
\bi
\i $\gw{\phi^t\#\mu,\phi^t\#\nu}\leq e^{Lt}\gw{\mu,\nu}$,
\i $\gw{\mu,\phi^t\#\mu}\leq b t\|v\|_{C^0} |\mu|$,
\i $\gw{\phi^t\#\mu,\psi^t\#\nu}\leq e^{Lt}\gw{\mu,\nu}+\frac{e^{Lt}-1}{L}\,|\mu|\,\sup_{\tau\in[0,t]}\|v_t-w_t\|_{C^0}$.
\ei
\ep
\bproof We first prove properties for the standard Wasserstein distance.

\b{Property 1}.  Let $\pi$ be the transference plan realizing $W_p\Pt{\mu,\nu}$. Observe that $\phi^t$ is a diffeomorphism of the space $\R^d$, then $\phi^t\times\phi^t$ is a diffeomorphism of the space $\R^d\times \R^d$. Since $\pi$ is a probability density on $\R^d\times\R^d$, then one can define $\pi':=(\phi^t\times\phi^t)\#\pi$, another probability density on $\R^d\times\R^d$. It is easy to prove that $\pi'$ is indeed a transference plan between $\phi^t\#\mu$ and $\phi^t\#\nu$. Then we can use such transference plan $\pi'$ to estimate $W_p\Pt{\phi^t\#\mu,\phi^t\#\nu}$. This gives
\bqn
W_p^p\Pt{\phi^t\#\mu,\phi^t\#\nu}&\leq& |\mu|^p\int |x-y|^p\,d\pi'(x,y)=|\mu|^p\int |\phi^t(x)-\phi^t(y)|^p\,d\pi(x,y)\leq\nn
&\leq& |\mu|^p\int e^{Lpt} |x-y|^p \,d\pi(x,y)=e^{Lpt} W^p_p(\mu,\nu),
\eqnn
where we used the definition of the push-forward in the first equality and the Gronwall lemma in the last inequality.

\b{Property 2}. Define the transference plan $\pi$ such that $d\pi(x,y)=|\mu|^{-1}d\mu(x)\delta_{y=\phi^t(x)}$ on $\R^d\times \R^d$. Observe that it is a transference plan between $\mu$ and $\phi^t\#\mu$. Then we have
\bqn
W_p^p\Pt{\mu,\phi^t\#\mu}&\leq& |\mu|^p\int |x-y|^p \,d\pi(x,y)=|\mu|^p\int |x-\phi^t(x)|^p \,|\mu|^{-1}\,d\mu(x)\leq |\mu|^p \int(\|v\|_{C^0}t)^p \,|\mu|^{-1}\,d\mu(x)=\nn
&=&|\mu|^p(\|v\|_{C^0}t)^p.
\eqnn

\b{Property 3}.  The proof is similar to proof of Property 1. Let $\pi$ be the transference plan realizing $W_p\Pt{\mu,\nu}$. Observe that $\phi^t\times\psi^t$ is a diffeomorphism of the space $\R^d\times \R^d$. Since $\pi$ is a probability density on $\R^d\times\R^d$, then one can define $\pi':=(\phi^t\times\psi^t)\#\pi$, another probability density on $\R^d\times\R^d$. It is easy to prove that $\pi'$ is indeed a transference plan between $\phi^t\#\mu$ and $\psi^t\nu$. Then we can use such transference plan $\pi'$ to estimate $W_p\Pt{\phi^t\#\mu,\psi^t\#\nu}$. We have
\bqn
W_p^p\Pt{\phi^t\#\mu,\psi^t\#\nu}&\leq& |\mu|^p\int |x-y|^p\,d\pi'(x,y)=|\mu|^p\int |\phi^t(x)-\psi^t(y)|^p\,d\pi(x,y)\leq\nn
&\leq& |\mu|^p\int \Pt{e^{Lt} |x-y|+\frac{e^{Lt}-1}{L} \sup_{\tau\in[0,t]}\|v_t-w_t\|_{C^0}}^p \,d\pi(x,y),
\eqnn
where we have used Gronwall inequality. Minkowski inequality now gives
\bqn
W_p\Pt{\phi^t\#\mu,\psi^t\#\nu}&\leq& |\mu| e^{Lt} \Pt{\int  |x-y|^p\,d\pi(x,y)}^{1/p}+|\mu|\frac{e^{Lt}-1}{L} \sup_{\tau\in[0,t]}\|v_t-w_t\|_{C^0}\,\int \,d\pi(x,y)=\nn
&=&e^{Lt} W_p(\mu,\nu)+|\mu| \frac{e^{Lt}-1}{L} \sup_{\tau\in[0,t]}\|v_t-w_t\|_{C^0},
\eqnn
where we also used $\int \,d\pi(x,y)=1$.

We now prove equivalent properties for the generalized Wasserstein distance.

\b{Property 1}.  Let $\tilde\mu\sotto\mu,\tilde\nu\sotto\nu$ be the choices realizing $\tw{\mu,\nu}$, i.e.
$$\tw{\mu,\nu}=a^p(|\mu-\tilde\mu|+|\nu-\tilde\nu|)^p+b^pW^p_p(\tilde\mu,\tilde\nu).$$ Then estimate $\tw{\phi^t\#\mu,\phi^t\#\nu}$ with  $\phi^t\#\tilde\mu$ and $\phi^t\#\tilde\nu$. Observe that $\phi^t\#\tilde\mu\sotto\phi^t\#\mu,\phi^t\#\tilde\nu\sotto\phi^t\#\nu$, and in particular $|\phi^t\#\mu-\phi^t\#\tilde\mu|=|\mu-\tilde\mu|$, and similarly for the other term. We then have
\bqn
\tw{\phi^t\#\mu,\phi^t\#\nu}&\leq& a^p(|\phi^t\#\mu-\phi^t\#\tilde\mu|+|\phi^t\#\nu-\phi^t\#\tilde\nu|)^p+b^pW^p_p(\phi^t\#\tilde\mu,\phi^t\#\tilde\nu)\leq\nn
&\leq&a^p (|\mu-\tilde\mu|+|\nu-\tilde\nu|)^p+ b^p e^{Lpt}W_p^p(\tilde\mu,\tilde\nu)\leq\nn
&\leq&e^{Lpt}\Pt{a^p (|\mu-\tilde\mu|+|\nu-\tilde\nu|)^p+ b^p W_p^p(\tilde\mu,\tilde\nu)}.
\eqnn
Computing the $p$-th root, we have the result. Proof of \b{Property 3} is completely equivalent, by using $\psi^t\#\tilde\nu\sotto\psi^t\#\nu$ and the corresponding inequality for $W_p(\phi^t\#\tilde\mu,\psi^t\#\tilde\nu)$.

\b{Property 2}. To estimate $\gw{\mu,\phi^t\#\mu}$, choose $\tilde\mu=\mu,\tilde\nu=\phi^t\#\mu$. Then one has $\gw{\mu,\phi^t\#\mu}\leq bW_p(\mu,\phi^t\#\mu)$. Using Property 2 for the standard Wasserstein distance, one has the result.
\eproof

\section{The generalized Wasserstein distance $W_1^{1,1}$ is the flat metric}
\label{s-flat}
\renewcommand{\gwbase}{W^{1,1}_1}
\renewcommand{\gw}[1]{W^{1,1}_{1}\Pt{#1}}

In this section, we provide a dual formulation for the generalized Wasserstein ditance $W_1^{1,1}$, proving that it coincides with the flat metric.
First define the spaces $\LL,\KK,\FF$ as follows:
$$\LL:=\Pg{f\in C^0_c(\R^d,\R)\ |\ \|f\|_\infty\leq 1},\qquad
\KK:=\Pg{f\in C^0_c(\R^d,\R)\ |\ \|f\|_{Lip}\leq 1},\qquad
\FF:=\LL\cap\KK.$$

We also recall the following dual formulation for $L^1$ and $W_1$ distances.
\bp \label{p-duali}
For all $\mu,\nu\in\Mu$ it holds $$|\mu-\nu|=\sup\Pg{\int f d(\mu-\nu)\ |\ f\in \LL}.$$

For all $\mu,\nu\in\Mu$ with $|\mu|=|\nu|$ it holds $$W_1(\mu,\nu)=\sup\Pg{\int f d(\mu-\nu)\ |\ f\in \KK}.$$
\ep
The second statement of Proposition \ref{p-duali} is known as the Kantorovich-Rubinstein theorem, see \cite[Theorem 1.14]{villani}.

We now recall the definition of the flat metric.
\bdeff
Let $\mu,\nu\in \Mu$. Define $$d(\mu,\nu):=\sup\Pg{\int f d(\mu-\nu)\ |\ f\in \FF}.$$
The functional $d$ is a metric on $\Mu$, called the \b{flat metric}.
\edeff

We now state the main result of this section.
\bt \label{t-flat}
Let $\mu,\nu\in\Mu$. Then 
\bqn\gw{\mu,\nu}=d(\mu,\nu).\eqnl{e-flat}
\et
The proof is based on some duality properties of convex functionals. For this reason, we first recall some useful definitions and results. For a complete description, see e.g. \cite{rocka}. In particular Theorem \ref{t-rocka} is Theorem 20.e. in \cite{rocka}.

\bdeff Let $X$ be a Banach space and $F:X\to \bar{\R}$ a function. The conjugate function $F^*:X^*\to\bar{\R}$ is 
$$F^*(y):=\sup_{x\in X}(\Pa{y,x}-f(x)).$$
\edeff
\bt\label{t-rocka} Let $X$ be a Banach space. Let $F_1,F_2:X\to\R\cup\Pg{+\infty}$ be convex and closed. Assume that there exists a neighborhood $U$ of the origin in $X$, an open set $M$ in $X^*$ and a constant $k$ such that for all sets 
$$V_\alpha:=\Pg{(y_1,y_2)\ |\ y_i\in \mathrm{dom}(F_i^*),\ y_1+y_2\in M,\ F_1^*(y_1)+F_2^*(y_2)\leq \alpha}$$
it holds 
\bqn
\sup_{x\in U, (y_1,y_2)\in V_\alpha}\Pa{y_1+y_2,x} <k.
\eqnl{e-crock}

Then the conjugate function $F^*$ of $F=F_1+F_2$ satisfies
\bqn
F^*(y)=\min_{y_1+y_2=y}\Pt{F_1^*(y_1)+F_2^*(y_2)}.
\eqnl{e-crock2}
\et
We recall that a function is closed if the set $\Pg{f\leq k}$ is closed for all $k\in\bar{R}$. Also observe that we removed $-\infty$ from the codomain of $F_1,F_2$. This gives that $F_1,F_2$ are both proper in the sense of \cite[p. 1]{rocka}.

\begin{proof}[Proof of Theorem \ref{t-flat}.] We define the following functionals on $X^*=(C^0_c(\R^n),\|\cdot\|_{\infty})$:
$$F_1(f):=
\begin{cases}
0&\mbox{~~~when~}\|f\|_\infty\leq 1\\
+\infty&\mbox{~~~elsewhere.}
\end{cases}
\qquad\mbox{and}\qquad
F_2(f):=
\begin{cases}
0&\mbox{~~~when~}\|f\|_{Lip}\leq 1\\
+\infty&\mbox{~~~elsewhere.}
\end{cases}
$$
Recall that the dual space $X^*$ is the space of signed Radon measures, see e.g. \cite[p.49]{ev-gar}. Then, dual formulations in Proposition \ref{p-duali} easily give that $F_1^*(\mu-\nu)=|\mu-\nu|$ and $F_2^*(\mu-\nu)=W_1(\mu,\nu)$. We now consider $F=F_1+F_2$ and study $F^*(\mu-\nu)$: it is easy to prove that it coincides with $d(\mu,\nu)$, by the definition of the conjugate function.

We now prove that $F^*(\mu-\nu)$ coincides with $W_1^{1,1}(\mu,\nu)$, by using Theorem \ref{t-rocka}. It is easy to prove that $F_1,F_2$ are proper, closed and convex functions, and that $U:=\Pg{\|f\|\mu\|_{\infty}<\eps}$, $M=\Pg{|\mu|< \eps}$, $k=\eps^2$ satisfy \r{e-crock}. Then, condition \r{e-crock2} reads as $F^*(\mu-\nu)=\min_{(\mu_1-\nu_1)+(\mu_2-\nu_2)=\mu-\nu}\Pt{|\mu_1-\nu_1|+W_1(\mu_2,\nu_2}$, that clearly coincides with $W_1^{1,1}(\mu,\nu)$.
\eproof

\section{Generalized Benamou-Brenier formula}
\label{s-bb}

\renewcommand{\gwbase}{W^{a,b}_2}
\renewcommandx{\gw}[1]{W^{a,b}_{2}\Pt{#1}}
\renewcommand{\twbase}{T^{a,b}_2}
\renewcommandx{\tw}[1]{T^{a,b}_{2}\Pt{#1}}
\renewcommand{\DD}{\mathcal{D}}
\renewcommand{\LL}{\mathcal{L}}
\newcommand{\RR}{\mathcal{R}}

In this section we generalize the Benamou-Brenier formula (recalled below, see \cite{ben-bren}) to $\gwbase$. The interest of such formula is to relate the Wasserstein distance between two measures $\mu_0,\mu_1$ to the minimization of the functional $\int |v_t|^2 d\mu_t$ among all solutions of the linear transport equation from $\mu_0$ to $\mu_1$. We first recall the original Benamou-Brenier formula. Observe that we deal with probability measures in $\Muu$.
\bt \label{t-bb} Let $\mu_0,\mu_1\in\Pac$ where $\Pac:=\Muu\cap \P$ is the space of probability measures that are absolutely continuous with respect with the Lebesgue measure and with compact support. Endow $\Pac$ with the weak-$*$ topology.

Let $V(\mu_0,\mu_1)$ be the set of couples measure-velocity field $(\mu,v):=(\mu_t,v_t)_{t\in\Pq{0,1}}$ such that $\mu\in C(\Pq{0,1},\Pac)$, $v\in L^2(d\mu_t dt)$, $\cup_{t\in\Pq{0,1}} \supp(\mu_t)$ is bounded, and such that they satisfy the following boundary value problem
$$\begin{cases}
\partial_t\mu_t+\nabla\cdot(v_t \mu_t)=0\\
\mu_{|_{t=0}}=\mu_0,\qquad \mu_{|_{t=1}}=\mu_1.
\end{cases}$$

Define the action functional $\fa{\mu,v}:=\int_0^1 dt\Pt{ \int_{\R^d} d\mu_t\,|v_t|^2}$ on $V(\mu_0,\mu_1)$. Then, it holds
\bqn W_2^2(\mu_0,\mu_1)=\inf\Pg{\fa{\mu,v}\,|\, (\mu,v)\in V(\mu_0,\mu_1)}.
\eqnl{e-bb}
\et
Such result has been proven to hold also in the larger space of probability measures with finite second order moments, see \cite{gradient}. It is also easy to prove that \r{e-bb} holds for $\mu_0,\mu_1\in \Muu$ with the same mass $m$. Indeed, it is sufficient to use \r{e-bb} for $m^{-1}\mu_0,m^{-1}\mu_1$ and to observe that we have the same degree of homogeneity on the left and right hand sides when multiplying by a constant.

We now prove that a similar result holds for $\gwbbase$ and the transport equation with source. We first define the space and the functional that we study.

\bdeff \label{d-a} Consider $\mu_0,\mu_1\in\Muu$. Let $V(\mu_0,\mu_1)$ be the set of triples (measure, velocity field, source term) $(\mu,v,h):=(\mu_t,v_t,h_t)_{t\in\Pq{0,1}}$ with the 
following properties:
$\mu\in C(\Pq{0,1},\Muu)$, with $\Muu$ endowed with the weak-$*$ topology;
$v\in L^2(d\mu_t dt)$;
$h\in L^1(\Pq{0,1},\Muu)$ in the sense that $\int_0^1 dt\Pt{\int_{\R^d} d|h_t|}<\infty$;
$\cup_{t\in\Pq{0,1}} \supp(\mu_t)$ is bounded;
they satisfy the following boundary value problem:
\bqn
\begin{cases}
\partial_t\mu_t+\nabla\cdot(v_t \mu_t)=h_t,\\
\mu_{|_{t=0}}=\mu_0,\qquad \mu_{|_{t=1}}=\mu_1.
\end{cases}
\eqnl{e-bvp}
We define the action functional on $V(\mu_0,\mu_1)$ by
$$\fb{\mu,v,h}:=a^2\Pt{\int_0^1 dt \Pt{\int_{\R^d}d |h_t|}}^2 +b^2\int_0^1 dt \Pt{ \int_{\R^d} d\mu_t\,|v_t|^2}.$$
\edeff

\brem \label{r-h} Observe that the conditions given above also imply that $\cup_{t\in[0,1]}\mathrm{supp}(h_t)\subset \cup_{t\in\Pq{0,1}} \supp(\mu_t)$, and in particular $h_t$ have uniformly bounded support. Indeed, by contradiction, assume that $\cup_{t\in[0,1]}\mathrm{supp}(h_t)\not\subset \cup_{t\in\Pq{0,1}} \supp(\mu_t)$. Looking at $h$ as a functional on $C^\infty_0$ functions, this means that there exists a function $\psi\in C^\infty_0\Pt{[0,1]\times\R^d,\R}$ with $\supp(\psi)\subset [0,1]\times \Pt{\R^d\setminus \Pt{\cup_{t\in\Pq{0,1}} \supp(\mu_t)}}$ and such that $\int_0^1 dt \int_{\R^d} dh_t \psi(t,x)\neq 0$. Observe now that, by construction, one has $\int_0^1 dt \int d \mu_t (\partial_t \psi+ v\cdot \nabla \psi)=0$, since $\psi$ and its derivatives are identically 0 on the support of $\mu_t$ for each $t\in [0,1]$. Observe now that $(\mu,v,h)$ satisfy \r{pde} in the weak sense. Choosing $\psi$ as a test function, one has $0=\int_0^1 dt \int_{\R^d} dh_t \psi(t,x)$. Contradiction.
\erem

We now state the generalized Benamou-Brenier formula:
\bt \label{t-gbb}
Let $\mu_0,\mu_1\in\Muu$. Then 
\bqn
\inf\Pg{\fb{\mu,v,h}\ |\ (\mu,v,h)\in V(\mu_0,\mu_1)}=\tw{\mu_0,\mu_1}.
\eqnl{e-gbb}
\et
It is clear the similarity between $\fbbase$ and $\fabase$. In particular, the standard Benamou-Brenier formula can be recovered as a particular case of Theorem \ref{t-gbb} when $h\equiv 0$ and $a\to \infty$.

\brem \label{r-altrib} It is possible to find a result similar to Theorem \ref{t-gbb} by changing the definition of both $\twbase$ and $\fbbase$. In particular, one can replace $\twbase$ with
$$\inf_{\tilde\mu,\tilde\nu\in\Mu,\,|\tilde\mu|=|\tilde\nu|}a^2\Pt{|\mu-\tilde\mu|^2+|\nu-\tilde\nu|^2}+b^2W_2^2(\tilde\mu,\tilde\nu),$$
and $\fbbase$ with
$$a^2\Pt{\int_0^1 dt \Pt{\int_{\R^d}d h_t^+}}^2 +a^2\Pt{\int_0^1 dt \Pt{\int_{\R^d}d h_t^-}}^2 + b^2\int_0^1 dt \Pt{ \int_{\R^d} d\mu_t\,|v_t|^2}.$$
This means that we have distributed the power 2 on the terms for creation and removal of mass, both for $\twbase$ and $\fbbase$. Proofs given below for Theorem \ref{t-gbb} can be easily adapted to this setting.
\erem

\bproof[Proof of Theorem \ref{t-gbb}] The proof is divided in 4 steps.

\b{Step 1.} We first prove the inequality $\fb{\mu,v,h}\geq \tw{\mu_0,\mu_1}$ under the following stronger regularity assumptions for $v,h$:
\bi
\i $v$ is uniformly $L$-Lipschitz with respect to $x$; it has $C^0$-norm uniformly bounded in time, i.e. $M:=\sup_{t\in\Pq{0,1}}\|v_t\|_{C^0}<\infty$;
\i $h\in L^\infty([0,1],\Muu)$, i.e. it satisfies $P:=\sup_{t\in\Pq{0,1}}\int_{\R^d}d|h_t(.)|<\infty$.
\ei  The idea of the proof is to approximate solutions of \r{pde} via an adapted sample-and-hold method, and to prove the inequality $\fb{\mu,v,h}\geq \tw{\mu_0,\mu_1}$ for such approximations.

The proof is divided into two substeps. Before the main parts of the proof, we state some simple remarks. First of all, since we deal with approximations of the dynamics given by $v,h$, then the approximated solution $\bmu$ could fail to be a positive measure for some times. Then, one needs to replace $\bmu$ with its positive part all along the proof. For simplicity of notation, this replacement is implict all along the proof.

Second, we fix some notations that will be useful all along the proof. Given the initial datum $\mu_0$, we will prove that all measures studied in the proof have bounded mass, and in particular $|\mu_t|,|\bmu_t|,|\bbmu_t|\leq |\mu_0|+P$. We define $$m:=|\mu_0|+P.$$
We also define
$$\al:=\sqrt2 \max\Pg{\frac{M}a,\frac{L}{b}},\qquad\beta:=2aP+bMm.$$

\b{Step 1.1:} In this step, we define an approximate solution $\bmu$, together with $v^k,h^k$, via a sample and hold method. We will prove that both  $\bmu\to\mu$ and $\fb{\bmu,v^k,h^k}\to \fb{\mu,v,h}$ for $k\to\infty$. 

Fix $k\in\N$ and define $\dt:=2^{-k}$. We discretize the time interval $[0,1]$ in small intervals $[n\dt,(n+1)\dt]$. The idea of the discretization is first to divide each interval $\Pq{n\dt,(n+1)\dt}$ in three parts: $$\Pq{n\dt, n\dt +\dt^2},\ \Pq{ n\dt +\dt^2,(n+1)\dt-\dt^2},\ \Pq{(n+1)\dt-\dt^2,(n+1)\dt}.$$ On the first part we use the negative part $h^-$ of $h$, then the velocity $v$, then the positive part $h^+$ of $h$. Clearly, each term must be correctly rescaled, to have $\bmu_{(n+1)\dt}$ close to $\mu_{(n+1)\dt}$.

We define the following vector field and the source term:
\bqn
v^k_{n\dt+\tau}:=\begin{cases}
\frac{\dt}{\dt-2\dt^2}v_{n\dt+\frac{\dt}{\dt-2\dt^2}(\tau-\dt^2)}&\mbox{~~for~~}\tau\in(\dt^2,\dt-\dt^2],\\
0 &\mbox{~~for~~}\tau\in(0,\dt^2]\cup(\dt-\dt^2,\dt],
\end{cases}
\eqnn
\bqn
h^k_{n\dt+\tau}:=\begin{cases}
\dt^{-1} h^-_{n\dt+\dt^{-1}\tau}&\mbox{~~for~~}\tau\in(0,\dt^2],\\
0&\mbox{~~for~~}\tau\in(\dt^2,\dt-\dt^2],\\
\dt^{-1} h^+_{n\dt+\dt^{-1}(\tau-(\dt-\dt^2))} &\mbox{~~for~~}\tau\in(\dt-\dt^2,\dt].
\end{cases}
\eqnn

 Observe that $v^k$ and $h^k$ will never act at the same time, i.e. $v^k_t\neq 0$ implies $h^k=0$ and viceversa. A scheme of the evolution of the mass $|\bmu_t|$ is given in Figure \ref{f-bmu}.

\newcommand{\spez}[4]{
\pgfxyline(#1)(#2)
\pgfxyline(#2)(#3)
\pgfxyline(#3)(#4)
}

\begin{figure}
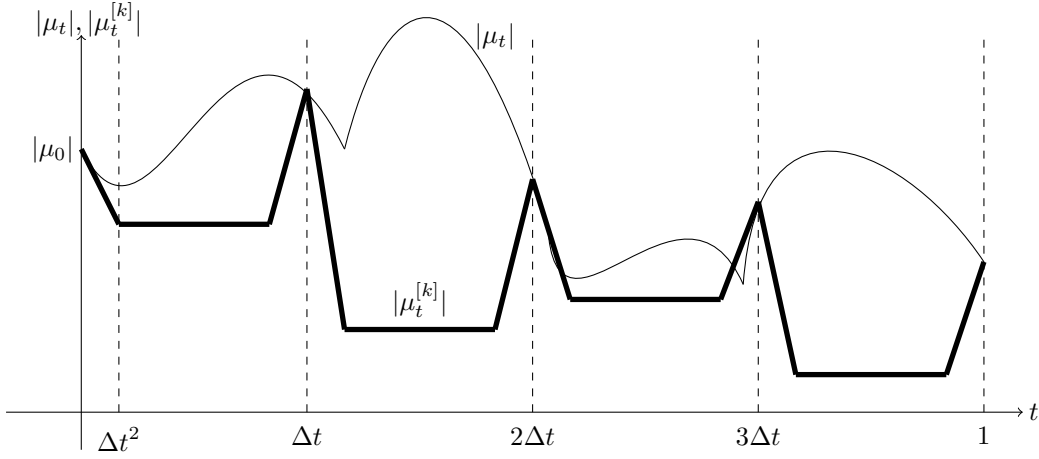


\begin{pgfpicture}{0cm}{0cm}{15cm}{6cm}


\begin{pgfscope}
\pgfsetendarrow{\pgfarrowto}
\pgfline{\pgfxy(0,.5)}{\pgfxy(13.5,.5)}
\pgfputat{\pgfxy(13.6,.5)}{\pgfbox[left,center]{$t$}}
\pgfputat{\pgfxy(1.5,0.3)}{\pgfbox[center,top]{$\dt^2$}}
\pgfputat{\pgfxy(4,0.3)}{\pgfbox[center,top]{$\dt$}}
\pgfputat{\pgfxy(7,0.3)}{\pgfbox[center,top]{$2\dt$}}
\pgfputat{\pgfxy(10,0.3)}{\pgfbox[center,top]{$3\dt$}}
\pgfputat{\pgfxy(13,0.3)}{\pgfbox[center,top]{$1$}}
\pgfline{\pgfxy(1,0)}{\pgfxy(1,5.5)}
\pgfputat{\pgfxy(1.1,5.5)}{\pgfbox[center,bottom]{$|\mu_t|,|\bmu_t|$}}
\pgfputat{\pgfxy(.9,4)}{\pgfbox[right,center]{$|\mu_0|$}}
\end{pgfscope}

\begin{pgfscope}
\pgfsetdash{{3pt}{3pt}}{0pt}
\pgfxyline(1.5,0.5)(1.5,5.5)
\pgfxyline(4,0.5)(4,5.5)
\pgfxyline(7,0.5)(7,5.5)
\pgfxyline(10,0.5)(10,5.5)
\pgfxyline(13,0.5)(13,5.5)
\end{pgfscope}

\pgfxycurve(1,4)(2,2)(3,7)(4.5,4)
\pgfxycurve(4.5,4)(5,6)(6,7)(7.2,3)
\pgfxycurve(7.2,3)(7.3,1)(9,4)(9.8,2.2)
\pgfxycurve(9.8,2.2)(10,5)(12,4)(13,2.5)
\pgfputat{\pgfxy(6.5,5.5)}{\pgfbox[center,center]{$|\mu_t|$}}

\begin{pgfscope}
\pgfsetlinewidth{2pt}
\spez{1,4}{1.5,3}{3.5,3}{4,4.8}
\spez{4,4.8}{4.5,1.6}{6.5,1.6}{7,3.6}
\spez{7,3.6}{7.5,2}{9.5,2}{10,3.3}
\spez{10,3.3}{10.5,1}{12.5,1}{13,2.5}
\pgfputat{\pgfxy(5.5,2)}{\pgfbox[center,center]{$|\bmu_t|$}}
\end{pgfscope}

\end{pgfpicture}

\caption{Evolution of $|\mu_t|,|\bmu_t|$ for $k=2$.}
\label{f-bmu}
\end{figure}

We now define $\bmu$ as the solution of \r{pde} in $C([0,1],\Muu)$ with velocity field $v^k$, source $h^k$, and initial datum $\bmu_0=\mu_0$.  It is evident that the measure has uniformly bounded mass, in particular $|\bmu_t|\leq m$ for all $t\in[0,1]$.

It is also easy to prove the following property: for $\tau\in[0,\dt]$ it holds
\bqn
\gw{\bmu_{n\dt},\bmu_{n\dt+\tau}}\leq\dt(2aP+bMm)=:\beta \dt.
\eqnl{e-quasilip}

We now prove that $\bmu$ is a Cauchy sequence with respect to the distance $\wwbase$ defined as follows
$$ \ww{\mu,\nu}=\sup_{t\in[0,1]}\gw{\mu_t,\nu_t}.$$
We recall that $C([0,1],\Mu)$ is complete with respect to $\wwbase$, as a direct consequence of the completeness of $\Mu$ with respect to $\gwbase$, see Proposition \ref{p-complete}.

First observe that, by substitution, the following formula holds for $\bmu_{(n+2)\dt}$:
\bqn
\bmu_{(n+2)\dt}&=&\ff{(n+1)\dt,(n+2)\dt}\#\bmu_{(n+1)\dt}-\ff{(n+1)\dt,(n+2)\dt}\#\hhm{(n+1)\dt,(n+2)\dt}+\hhp{(n+1)\dt,(n+2)\dt}=\nn
&=&\ff{(n+1)\dt,(n+2)\dt}\#\Pt{\ff{n\dt,(n+1)\dt}\#\Pt{\bmu_{n\dt}-\hhm{n\dt,(n+1)\dt}}+\hhp{n\dt,(n+1)\dt}}+\nn
&&-\ff{(n+1)\dt,(n+2)\dt}\#\hhm{(n+1)\dt,(n+2)\dt}+\hhp{(n+1)\dt,(n+2)\dt}=\nn
&=&\ff{n\dt,(n+2)\dt}\#\bmu_{n\dt}-\ff{n\dt,(n+2)\dt}\#\hhm{n\dt,(n+1)\dt}+ \ff{(n+1)\dt,(n+2)\dt}\# \hhp{n\dt,(n+1)\dt}+\nn
&&-\ff{(n+1)\dt,(n+2)\dt}\#\hhm{(n+1)\dt,(n+2)\dt}+\hhp{(n+1)\dt,(n+2)\dt},
\eqnn
where 
\bi
\i $\ff{t_1,t_2}$ is the diffeomorphism corresponding to the flow generated by $v^k$ on the time interval $[t_1,t_2]$;
\i $\hhm{t_1,t_2}:=\int_{t_1}^{t_2} \underline{h}^k_t\,dt$ is the mass removal given by $\underline{h}^k$ on the time interval $[t_1,t_2]$;
\i $\hhp{t_1,t_2}:=\int_{t_1}^{t_2} \overline{h}^k_t\,dt$ is the mass creation given by $\overline{h}^k$ on the time interval $[t_1,t_2]$.
\ei
We also decompose $\bmu[k-1]_{(n+2) \dt}$ by using properties of composition of 
$\Phi^k,\underline{H}^k,\overline{H}^k$. This gives:
\bqn
\bmu[k-1]_{(n+2) \dt}&=&\ff{n\dt,(n+2)\dt}\#\bmu[k-1]_{n\dt}-\ff{n\dt,(n+2)\dt}\#\hhm{n\dt,(n+1)\dt}+\nn
&&-\ff{(n+1)\dt,(n+2)\dt}\#\ff{n\dt,(n+1)\dt}\#\hhm{(n+1)\dt,(n+2)\dt}+\hhp{n\dt,(n+1)\dt}+\hhp{(n+1)\dt,(n+2)\dt}
\eqnn

We now estimate $\gw{\bmu[k-1]_{(n+2) \dt},\bmu_{(n+2) \dt}}$ with respect to $\gw{\bmu[k-1]_{n\dt},\bmu_{n\dt}}$, i.e. the value of $\gwbase$ at the right extreme of the interval of discretization for $k-1$ with respect to its value at the left extreme. We choose $n$ even. Using estimates in Proposition \ref{p-flow}, we have:
\bqn
&&\gw{\bmu[k-1]_{(n+2) \dt},\bmu_{(n+2) \dt}}\leq \gw{\ff{n\dt,(n+2)\dt}\#\bmu[k-1]_{n\dt},\ff{n\dt,(n+2)\dt}\#\bmu_{n\dt}}+\nn
&&+\gw{\ff{n\dt,(n+2)\dt}\#\hhm{n\dt,(n+1)\dt},\ff{n\dt,(n+2)\dt}\#\hhm{n\dt,(n+1)\dt}}+\nn
&&+\gw{\ff{(n+1)\dt,(n+2)\dt}\#\ff{n\dt,(n+1)\dt}\#\hhm{(n+1)\dt,(n+2)\dt},\ff{(n+1)\dt,(n+2)\dt}\#\hhm{(n+1)\dt,(n+2)\dt}}+\nn
&&+\gw{\hhp{n\dt,(n+1)\dt},\ff{(n+1)\dt,(n+2)\dt}\# \hhp{n\dt,(n+1)\dt}}+\gw{\hhp{(n+1)\dt,(n+2)\dt},\hhp{(n+1)\dt,(n+2)\dt}}\leq\nn
&&\leq e^{L\dt}\gw{\bmu[k-1]_{n\dt},\bmu_{n\dt}}+0+e^{L\dt}\gw{\ff{n\dt,(n+1)\dt}\#\hhm{(n+1)\dt,(n+2)\dt},\hhm{(n+1)\dt,(n+2)\dt}}+\nn
&&+b\dt M(\dt P)+0\leq e^{ L\dt}\gw{\bmu[k-1]_{n\dt},\bmu_{n\dt}}+be^{ L\dt}\dt M  (\dt P)+b\dt M(\dt P).
\eqnl{e-chiave}

We apply the last inequality recursively. First recall that $\gw{\bmu[k-1]_0,\bmu_0}=0$ and that, for a sufficiently big $k$, it holds $e^{L\dt}\leq 1+2L\dt$ and $2L\dt\leq 1$. This gives 
\bqn
\gw{\bmu[k-1]_{n\dt},\bmu_{n\dt}}&\leq& b MP\dt^{2}(2+2L\dt)\frac{(1+2L\dt)^{n/2}-1}{1+2L\dt-1}\leq 2 bMP\dt \frac{e^{\frac nL\dt}-1}{2L}\leq\nn
&\leq& 2bMP2^{-k}\frac{e^{ L}-1}{L},
\eqnn
where we have used that $n\dt\leq 1$. Observe that the estimate is independent of $n$. Applying it recursively, one has
$$\gw{\bmu_{n\dt},\bmu[k+l]_{n\dt}}\leq 
\frac{2bMP(e^{ L}-1)}{L} 2^{-(k+1)} \frac{1-2^{-l/2}}{1-2^{-1/2}}.$$
Finally, take any $t\in\Pq{0,1}$: for each integer $k$, let $n_k$ be the biggest even number such that $n_k2^{-k}\leq t$. It clearly holds $|t-n2^{-k}|<2^{-k+1}$. One has
\bqn
&&\gw{\bmu_t,\bmu[k+l]_t} \leq \gw{\bmu_t,\bmu_{n_k2^{-k}}}+\gw{\bmu_{n_k2^{-k}},\bmu[k+l]_{n_k2^{-k}}}+
\gw{\bmu[k+l]_{n_k2^{-k}},\bmu[k+l]_t}\leq\nn
&&\leq 2\beta 2^{-k} +\frac{2bMP (e^{ L}-1)}{L}2^{-(k+1)}\frac{1-2^{-l/2}}{1-2^{-1/2}}+2\beta 2^{-k},
\eqnn
where we have used \r{e-quasilip} twice for the first term and $2^{l+1}$ times for third term. Since the estimate does not depend on $t$, one has
$d(\bmu,\bmu[k+l])\leq C_1 2^{-k}$ with $C_1:=4\beta+\frac{bMP (e^{L}-1)}{L} \frac{\sqrt{2}}{\sqrt{2}-1}$. Since the estimate does not depend on $l$ and $\ww{\bmu,\bmu[k+l]}\to 0$ for $k\to \infty$, we have that $\bmu$ is a Cauchy sequence. Since $C([0,1],\Mu)$ is complete with respect to $\wwbase$, then there exists a limit $\mu^*:=\lim_{k\to\infty}\bmu$, with $\mu^*\in\Mu$.

We now prove that $\mu^*=\mu$. We prove it by proving that it is a weak solution of \r{pde}. By uniqueness the result will follow. We have to prove that, for any\footnote{The index $t$ will be useful in the following change of variable in time.} $f_t\in C^\infty_0([0,1]\times \R^d)$, it holds
\bqn
\int_0^1 dt\Pt{\int_{\R^d}\Pt{d\mu^*_t(\partial_t f_t +v_t\cdot \nabla f_t)+dh_t f_t}}=0.
\eqnl{e-weak}
Observe that $\bmu$ is a solution of \r{pde} with vector field $v^k$, and source $h^k$. Then
\bqn
\int_0^1 dt\Pt{\int_{\R^d}\Pt{d\bmu_t(\partial_t f_t +v^k_t\cdot \nabla f_t)+dh^k_t f_t}}=0.
\eqnn
One can prove \r{e-weak} by proving the three following limits:
\vskip 5mm\noindent
1. $\lim_k \Pabs{\int_0^1\,dt\Pt{\int_{\R^d}\Pt{d\mu^*_t-d\bmu_t} \partial_t f_t}}=0$. This is a consequence of \r{e-gwintegrale}. Indeed, one has 
\bqn
\lim_k \Pabs{\int_0^1\, dt\Pt{\int_{\R^d}\Pt{d\mu^*_t-d\bmu_t} \partial_t f_t}}&\leq&
\int_0^1\,dt \Pt{\gw{\mu^*_t,\bmu_t}{\sqrt2}\max\Pg{\frac{\|\partial_tf_t\|_{\infty}}a,\frac{\|\partial_tf_t\|_{Lip}}b}}\leq\nn
&\leq& \ww{\mu^*,\bmu}{\sqrt2}\max\Pg{\frac{\|\partial_tf_t\|_{\infty}}a,\frac{\|\partial_tf_t\|_{Lip}}b}\to 0
\eqnn

\noindent
2. $\lim_k \Pabs{\int_0^1\,dt\Pt{\int_{\R^d}d\mu^*_t v_t \cdot \nabla f_t-d\bmu_t v^k_t \cdot \nabla f_t}}=0$. We first fix $k$ and $\dt:=2^{-k}$, and estimate
\bqn
\int_{0}^{\dt}\,dt\Pt{\int_{\R^d}d\mu^*_{n\dt+t} v_{n\dt+t} \cdot \nabla f_{n\dt+t}}-\int_{0}^{\dt}\,d\tau\Pt{\int_{\R^d}d\bmu_{n\dt+\tau} v^k_{n\dt+\tau} \cdot \nabla f_{n\dt+\tau}}.
\eqnl{e-p2}
Using the definition of $v^k_{n\dt+\tau}$, we have that it is 0 for $\tau\in[0,\dt^2]\cup (\dt-\dt^2,\dt]$ and that for $\tau\in(\dt^2,\dt-\dt^2]$ it holds $v^k_{n\dt+\tau}=\frac{\dt}{\dt-2\dt^2} v_{n\dt+\frac{\dt}{\dt-2\dt^2}(\tau-\dt^2)}$. Then, after the change of variable $\tau\to t:=(\tau-\dt^2)\frac{\dt}{\dt-2\dt^2}$, we have
\bqn
&&\int_{0}^{\dt}\,d\tau\Pt{\int_{\R^d}d\bmu_{n\dt+\tau} v^k_{n\dt+\tau} \cdot \nabla f_{n\dt+\tau}}=\nn
&&=\frac{\dt}{\dt-2\dt^2} \int_{\dt^2}^{\dt-\dt^2}\,d\tau\Pt{\int_{\R^d}d\bmu_{n\dt+\tau} v_{n\dt+\frac{\dt}{\dt-2\dt^2}(\tau-\dt^2)} \cdot \nabla f_{n\dt+\tau}}=\nn
&&=
\frac{\dt}{\dt-2\dt^2} \int_{0}^{\dt}\,\frac{\dt-2\dt^2}{\dt} dt \Pt{\int_{\R^d}d\bmu_{n\dt+\dt^2+ \frac{\dt-2\dt^2}{\dt}t} v_{n\dt+t} \cdot \nabla f_{n\dt+\dt^2+ \frac{\dt-2\dt^2}{\dt}t}}.
\eqnn
To go back to \r{e-p2}, we estimate for each $t\in[0,\dt]$ the following quantity\footnote{Here we denote with $\|\nabla f_t\|_{Lip}$ the Lipschitz constant for $\nabla f_t$ with respect to all $t,x$-variables, even if for \r{e-p21} the Lipschitz constant in space is needed only.}:
\bqn
&&\Pabs{\int_{\R^d}d\mu^*_{n\dt+t} v_{n\dt+t} \cdot \nabla f_{n\dt+t}-\int_{\R^d}d\bmu_{n\dt+\dt^2+ \frac{\dt-2\dt^2}{\dt}t} v_{n\dt+t} \cdot \nabla f_{n\dt+\dt^2+ \frac{\dt-2\dt^2}{\dt}t}}\leq\nn
&&\Pabs{\int_{\R^d}d\mu^*_{n\dt+t} v_{n\dt+t} \cdot \nabla f_{n\dt+t}-\int_{\R^d}d\bmu_{n\dt+\dt^2+ \frac{\dt-2\dt^2}{\dt}t} v_{n\dt+t} \cdot  \nabla f_{n\dt+t}}+\nn
&&+\Pabs{\int_{\R^d}d\bmu_{n\dt+\dt^2+ \frac{\dt-2\dt^2}{\dt}t} v_{n\dt+t} \cdot  \nabla f_{n\dt+t}-\int_{\R^d}d\bmu_{n\dt+\dt^2+ \frac{\dt-2\dt^2}{\dt}t} v_{n\dt+t} \cdot \nabla f_{n\dt+\dt^2+ \frac{\dt-2\dt^2}{\dt}t}}\leq\nn
&&\gw{\mu^*_{n\dt+t},\bmu_{n\dt+\dt^2+ \frac{\dt-2\dt^2}{\dt}t}}{\sqrt2}\max\Pg{\frac{M\|\nabla f_{t}\|_\infty}{a},\frac{L\|\nabla f_{t}\|_{Lip}}{b}}+\label{e-p21}\\
&&+\Pabs{\bmu_{n\dt+\dt^2+ \frac{\dt-2\dt^2}{\dt}t}}\, M\, \|\nabla f_{n\dt+t}-\nabla f_{n\dt+\dt^2+ \frac{\dt-2\dt^2}{\dt}t}\|_\infty.
\eqnn
We estimate the first term of the right hand side of \r{e-p21} via 
$$\gw{\mu^*_{n\dt+t},\bmu_{n\dt+\dt^2+ \frac{\dt-2\dt^2}{\dt}t}}\leq
\gw{\mu^*_{n\dt+t},\bmu_{n\dt+t}}+\gw{\bmu_{n\dt+t},\bmu_{n\dt+\dt^2+ \frac{\dt-2\dt^2}{\dt}t}}.$$
We estimate $\gw{\bmu_{n\dt+t},\bmu_{n\dt+\dt^2+ \frac{\dt-2\dt^2}{\dt}t}}$ by studying three cases:\\
(a) $t\in[0,\dt^2]$: We observe that the evolution from $\bmu_{n\dt+t}$ to $\bmu_{n\dt+\dt^2}$ is given by removal of mass $\hhm{n\dt+t,n\dt+\dt^2}$, while the evolution from $\bmu_{n\dt+\dt^2}$ to $\bmu_{n\dt+\dt^2+ \frac{\dt-2\dt^2}{\dt}t}$ is given by the push-forward of the diffeomorphism $\ff{n\dt+\dt^2,n\dt+\dt^2+ \frac{\dt-2\dt^2}{\dt}t}$. We then have 
\bqn
\gw{\bmu_{n\dt+t},\bmu_{n\dt+\dt^2+ \frac{\dt-2\dt^2}{\dt}t}}&\leq& \gw{\bmu_{n\dt+t},\bmu_{n\dt+\dt^2}}  +\nn
+\gw{\bmu_{n\dt+\dt^2},\bmu_{n\dt+\dt^2+ \frac{\dt-2\dt^2}{\dt}t}} &\leq&|t-\dt^2|\dt^{-1}P+ b\frac{\dt-2\dt^2}{\dt}t\|v^k\|_{C^0} m=\nn
&=&|t-\dt^2|\dt^{-1}P+ bM mt.\eqnl{e-cambio1}
(b) $t\in(\dt^2,\dt-\dt^2]$: We observe that the evolution is given by the push-forward of the diffeomorphism $\ff{n\dt+t,n\dt+\dt^2+ \frac{\dt-2\dt^2}{\dt}t}$. We have 
\bqn
\gw{\bmu_{n\dt+t},\bmu_{n\dt+\dt^2+ \frac{\dt-2\dt^2}{\dt}t}}&\leq& b\Pabs{t-\Pt{\dt^2+ \frac{\dt-2\dt^2}{\dt}t}} \|v^k\|_{C^0} m\leq\nn
&\leq& b|2t\dt -\dt^2|\frac{\dt}{\dt-2\dt^2}Mm.
\eqnn
(c) $t\in[\dt-\dt^2,\dt]$: This is similar to case 1. We have
\bqn
\gw{\bmu_{n\dt+t},\bmu_{n\dt+\dt^2+ \frac{\dt-2\dt^2}{\dt}t}}&\leq&
 |t-\dt+\dt^2|\dt^{-1}P+ \nn
 &&+b\Pabs{\frac{\dt-2\dt^2}{\dt}t-(\dt-2\dt^2)}  Mm.
\eqnl{e-cambio2}
We estimate the second term of the right hand side of \r{e-p21} via\footnote{Here it is sufficient to use the Lipschitz constant in the time variable.} $\Pabs{\bmu_{n\dt+\dt^2+ \frac{\dt-2\dt^2}{\dt}t}}\leq m$ and
$$\|\nabla f_{n\dt+t}-\nabla f_{n\dt+\dt^2+ \frac{\dt-2\dt^2}{\dt}t}\|_\infty\leq
\|\nabla f_t\|_{Lip}\Pabs{t-\Pt{\dt^2+ \frac{\dt-2\dt^2}{\dt}t}}=\|\nabla f_t\|_{Lip}\Pabs{2t\dt-\dt^2}.$$
Observe that both terms of the right hand side of \r{e-p21} have a symmetry property: the value in $t$ coincides with the value in $\dt-t$.\\
Back to \r{e-p2} and, by using \r{e-p21} and the symmetry described above, we have
\bqn
&&\Pabs{\int_{0}^{\dt}\,dt\Pt{\int_{\R^d}d\mu^*_{n\dt+t} v_{n\dt+t} \cdot \nabla f_{n\dt+t}}-\int_{0}^{\dt}\,d\tau\Pt{\int_{\R^d}d\bmu_{n\dt+\tau} v^k_{n\dt+\tau} \cdot \nabla f_{n\dt+\tau}}}\leq\nn
&&\leq 2{\sqrt2}\max\Pg{\frac{M\|\nabla f_{t}\|_\infty}{a},\frac{L\|\nabla f_{t}\|_{Lip}}{b}} \int_0^\dt dt\,\gw{\mu^*_{n\dt+t},\bmu_{n\dt+t}} +\nn
&&+2{\sqrt2}\max\Pg{\frac{M\|\nabla f_{t}\|_\infty}{a},\frac{L\|\nabla f_{t}\|_{Lip}}{b}} 
\int_0^{\dt^2}\,dt\Pt{ |t-\dt^2|\dt^{-1}P+ bMmt} +\nn
&& +2{\sqrt2}\max\Pg{\frac{M\|\nabla f_{t}\|_\infty}{a},\frac{L\|\nabla f_{t}\|_{Lip}}{b}}
\int_{\dt^2}^{\dt/2} \,dt\,b|2t\dt -\dt^2|\frac{\dt}{\dt-2\dt^2}Mm+\nn
&&+2\int_{0}^{\dt/2}\,dt\, m M \|\nabla f_t\|_{Lip}\Pabs{2t\dt-\dt^2}\leq C \ww{\mu^*,\bmu}\dt+ C \dt^3/2 +\nn
&&+ C(\dt-2\dt^2)\dt^3/2+C \int_0^{\dt/2}\,dt\, |2t\dt-\dt^2|  +C \int_0^{\dt/2}\,dt\, |2t\dt-\dt^2| 
\eqnl{e-p22}
with $C=2{\sqrt2}\max\Pg{\frac{M\|\nabla f_{t}\|_\infty}{a},\frac{L\|\nabla f_{t}\|_{Lip}}{b},\|\nabla f_t\|_{Lip}}\cdot
\max\Pg{1,P,2bMm, Mm}$. The estimate holds for $k\geq 2$, for which it holds $\frac{\dt}{\dt-2\dt^2}\leq 2$. We simply estimate \r{e-p22} with 
$C \ww{\mu^*,\bmu}\dt +C \dt^3/2+ C \dt^4/2+C \dt^3/2+C \dt^3/2<C \ww{\mu^*,\bmu}\dt+3 C \dt^3$, by using $|2 t\dt-\dt^2|\leq \dt^2$.\\
Going back to our estimate, using \r{e-p2} on each interval $[n\dt,(n+1)\dt]$, we have
\bqn
&&\lim_k \Pabs{\int_0^1\,dt\Pt{\int_{\R^d}d\mu^*_t v_t \cdot \nabla f_t-d\bmu_t v^k_t \cdot \nabla f_t}}\leq
 \lim_k \sum_{n=0}^{2^k-1} \left|\int_{0}^{\dt}\,dt\Pt{\int_{\R^d}d\mu^*_{n\dt+t} v_{n\dt+t} \cdot \nabla f_{n\dt+t}}+\right.\nn
&&\left.  -\int_{0}^{\dt}\,d\tau\Pt{\int_{\R^d}d\bmu_{n\dt+\tau} v^k_{n\dt+\tau} \cdot \nabla f_{n\dt+\tau}}\right|\leq   \lim_k 2^k (
C \ww{\mu^*,\bmu}2^{-k} +3 C 2^{-3k})=\nn
&&=\lim_k C \ww{\mu^*,\bmu}=0.
\eqnn

\noindent
3. $\lim_k \Pabs{\int_0^1\,dt\Pt{\int_{\R^d}\,d(h_t-h^k_t) f_t}}=0$. We first fix $k$ and $\dt:=2^{-k}$, and using again estimates in Proposition \ref{p-flow}, we have
\bqn
&&\lim_k \Pabs{\int_0^1\,dt \Pt{\int_{\R^d}\,d(h_t-h^k_t) f_t}}\leq \lim_k 2\cdot 2^k P\|f_t\|_{Lip}(1-\dt)\frac{2^{-2k}}2=0
\eqnn

We have proved that $\mu^*$ is a solution of \r{pde}, with $\mu^*\in C([0,1], \Mu)$. Observe now that $\mu^*-\mu$ is a solution of \r{pde} with initial datum $0$, vector field $v_t$ and source $0$. Applying standard result of existence and uniqueness of solutions of \r{pde} with zero source in $C([0,1], \Mu)$, we have $\mu^*=\mu$. Since $\mu\in C([0,1], \Muu)$, then $\mu^*\in C([0,1], \Muu)$ too.

We now prove that $\fb{\bmu,v^k,h^k}\to \fb{\mu,v,h}$ for $k\to \infty$. For the velocity term, we decompose
\bqn
&&\Pabs{\int_0^\dt dt\Pt{\int_{\R^d} d\mu_{n\dt+t} |v_{n\dt+t}|^2} -\int_0^\dt d\tau\Pt{\int_{\R^d} d\bmu_{n\dt+\tau} |v^k_{n\dt+\tau}|^2}}\leq\label{e-p4}\\
&& \Pabs{\int_0^\dt dt\Pt{\int_{\R^d} d\mu_{n\dt+t} |v_{n\dt+t}|^2} -\int_0^\dt dt\Pt{\int_{\R^d} d\bmu_{n\dt+t} |v_{n\dt+t}|^2}}+\nn
&&+ \Pabs{\int_0^\dt dt\Pt{\int_{\R^d} d\bmu_{n\dt+t} |v_{n\dt+t}|^2}- \int_0^\dt d\tau\Pt{\int_{\R^d} d\bmu_{n\dt+\tau} |v^k_{n\dt+\tau}|^2}}
\eqnn
We can easily estimate the first term by $$\int_0^\dt dt\,{\sqrt2}\max\Pg{\frac{M^2}a,\frac{2LM}{b}} \gw{\mu_t, \bmu_t}\leq  2M\al \ww{\mu,\bmu}\dt.$$
For the second term, we apply the change of variable $\tau\to t=(\tau-\dt^2)\frac{\dt}{\dt-2\dt^2}$ and find
\bqn
&&\Pabs{\int_0^\dt dt\Pt{\int_{\R^d} d\bmu_{n\dt+t} |v_{n\dt+t}|^2}- \int_0^\dt dt \Pt{\int_{\R^d} d\bmu_{n\dt+\dt^2+ \frac{\dt-2\dt^2}{\dt} t} |v_{n\dt+t}|^2}}\leq\nn
&&\leq 2M\al \int_0^\dt dt \gw{\bmu_{n\dt+t},\bmu_{n\dt+\dt^2+ \frac{\dt-2\dt^2}{\dt} t}}\leq 4M\al\beta \dt^2.
\eqnn
Going back to \r{e-p4}, we estimate the right-hand side with 
$2M\al \,  \ww{\mu,\bmu}\dt+4M\al\beta\, \dt^2\leq K  \ww{\mu,\bmu} \dt + K \dt^2,$
where $K:=\max\Pg{2M\al,4M\al\beta}$.\\
For the source part, the definition of $h^k$ easily gives 
$$\Pt{\int_0^1 dt \Pt{ \int_{\R^d}d|h_t|}}^2-\Pt{\int_0^1 dt \Pt{\int_{\R^d}d|h^k_t|}}^2=0.$$

Summing up, we have
\bqn
\Pabs{\fb{\bmu,v^k,h^k}-\fb{\mu,v,h}}\leq b^2 K \Pt{\ww{\mu,\bmu}+2^{-k}},
\eqnn
that gives $\lim_k \fb{\bmu,v^k,h^k}=\fb{\mu,v,h}$.

\b{Step 2:} We now define a $\bbmu$, together with $\bbv,\bbh$, that satisfies the three following properties:
\be
\i $\bbmu$ drives $\mu_0$ to $\bmu_1$, i.e. $(\bbmu,\bbv,\bbh)\in V(\mu_0,\bmu_1)$;
\i it holds $\fb{\bbmu,\bbv,\bbh}\leq\fb{\bmu,v^k,h^k}$;
\i it holds $\tw{\mu_0,\bmu_1}\leq \fb{\bbmu,\bbv,\bbh}$.
\ee
The idea is that, for each interval $[n\dt,(n+1)\dt]$ we move all the decreasing of mass in $[n\dt,n\dt+\dt^2)$, all the transport in $[n\dt+\dt^2,(n+1)\dt-\dt^2)$ and all the increase of mass in $[(n+1)\dt-\dt^2,(n+1)\dt]$. We divide this step in three substeps. In the first, we define $\bbmu$. In the second, we prove the properties stated above. In the third, we prove the result $\tw{\mu_0,\mu_1}\leq\fb{\mu,v,h}$ with the stronger regularity assumptions on $v,h$ recalled in Step 1.

\b{Step 2.1:} We now define $\bbmu$. With this goal, we define three transformations of measures. The transformation induced on the mass is described in Figure \ref{f-rld}.

\begin{figure}

\begin{pgfpicture}{0cm}{0cm}{15cm}{6cm}


\begin{pgfscope}
\pgfsetendarrow{\pgfarrowto}
\pgfline{\pgfxy(0,.5)}{\pgfxy(13.5,.5)}

\pgfline{\pgfxy(1,0)}{\pgfxy(1,5.5)}
\pgfputat{\pgfxy(1.1,5.5)}{\pgfbox[center,bottom]{$|\mu_t|$}}
\end{pgfscope}

\begin{pgfscope}
\pgfsetdash{{3pt}{3pt}}{0pt}
\pgfxyline(4,0.5)(4,5.5)
\pgfxyline(9.5,0.5)(9.5,5.5)
\pgfxyline(10.5,0.5)(10.5,5.5)
\pgfxyline(13,0.5)(13,5.5)
\pgfputat{\pgfxy(13.6,.5)}{\pgfbox[left,center]{$t$}}
\pgfputat{\pgfxy(4,0.3)}{\pgfbox[center,top]{$\bar t_1$}}
\pgfputat{\pgfxy(9.5,0.3)}{\pgfbox[center,top]{$\bar t_2$}}
\pgfputat{\pgfxy(10.5,0.3)}{\pgfbox[center,top]{$\bar t_3$}}
\pgfputat{\pgfxy(13,0.3)}{\pgfbox[center,top]{$1$}}
\end{pgfscope}

\spez{1,5}{1.5,4}{3.5,4}{4,4.8}
\spez{4,4.8}{4.5,3.6}{6.5,3.6}{7,4.6}
\spez{7,4.6}{7.5,3}{9.5,3}{10,2}
\spez{10,2}{10.5,2.5}{12.5,2.5}{13,3.5}

\begin{pgfscope}
\pgfsetlinewidth{2pt}
\spez{3.5,4}{4,2.8}{4,2.8}{4.5,3.6}
\pgfputat{\pgfxy(3.9,3)}{\pgfbox[right,center]{$\DD_{\bar t_1}$}}

\spez{7.5,3}{8,2}{8,2}{10,2}
\pgfputat{\pgfxy(9,1.8)}{\pgfbox[center,top]{$\LL_{\bar t_2}$}}

\spez{10,2}{12,2}{12,2}{12.5,2.5}
\pgfputat{\pgfxy(11.5,1.8)}{\pgfbox[center,top]{$\RR_{\bar t_3}$}}

\end{pgfscope}

\end{pgfpicture}

\caption{Transformations DOWN $\DD_{\bar t_1}$, LEFT $\LL_{\bar t_2}$ and RIGHT $\RR_{\bar t_3}$.}
\label{f-rld}
\end{figure}

\newcommand{\fff}[1]{\Phi_{[#1]}}

\noindent
\b{Transformation DOWN $\DD$}: The idea is to replace the increase-decrease of mass with the decrease-increase. Let $(\mu, v,h)$ be given, and $\bar t$ be a time such that: $ v_t=0$ on the interval $[\bar t-\dt^2,\bar t +\dt^2]$; $\bar h^-_t=0$ on the interval $[\bar t-\dt^2,\bar t]$; $\bar h^+_t=0$ on the interval $[\bar t,\bar t+\dt^2]$. Then replace $ h$ with $\hat h$ defined as follows:
\bqn
\hat h_t:=\begin{cases}
 h_t&\mbox{~~for $t\in[0,\bar t-\dt^2]\cup(\bar t +\dt^2,1]$},\\
 h_{t+\dt^2}&\mbox{~~for $t\in(\bar t-\dt^2,\bar t]$},\\
 h_{t-\dt^2}&\mbox{~~for $t\in(\bar t,\bar t +\dt^2]$}.
\end{cases}
\eqnn
Keep $ v$. We use the notation $\DD_{\bar t}$ for the solution $\hat \mu$ of \r{pde} with $v$ and $\hat h$, i.e. $\DD_{\bar t}( \mu):=\hat \mu$. We also denote $\DD_{\bar t}(\mu, v, h):=(\hat \mu,\hat v,\hat h)$.\\
\b{Transformation LEFT $\LL$}: The idea is to replace the transport-decrease with the decrease-transport. Let $(\mu, v, h)$ be given, and $\bar t$ be a time such that: $ h^+_t=0$ on the interval $[\bar t-\dt+2\dt^2,\bar t +\dt^2]$; $ h^-_t=0$ on the interval $[\bar t-\dt+2\dt^2,\bar t]$; $ v_t=0$ on the interval $[\bar t,\bar t+\dt^2]$. Then replace $ v$ with $\hat v$ defined as follows:
\bqn
\hat v_t:=\begin{cases}
 v_t&\mbox{~~for $t\in[0,\bar t-\dt+2\dt^2]\cup(\bar t +\dt^2,1]$},\\
0&\mbox{~~for $t\in(\bar t-\dt+2\dt^2,\bar t-\dt+3\dt^2]$},\\
 v_{t-\dt^2}&\mbox{~~for $t\in(\bar t-\dt+3\dt^2,\bar t +\dt^2]$}.
\end{cases}
\eqnn
Also replace $ h^-$ with $\hat h^-$ defined as follows:
\bqn
\hat h^-_t:=\begin{cases}
 h^-_t&\mbox{~~for $t\in[0,\bar t-\dt+2\dt^2]\cup(\bar t +\dt^2,1]$},\\
\Pt{\fff{\bar t-\dt+2\dt^2,\bar t}}^{-1}\# h^-_{t+\dt-2\dt^2}&\mbox{~~for $t\in(\bar t-\dt+2\dt^2,\bar t-\dt+3\dt^2]$},\\
0&\mbox{~~for $t\in(\bar t-\dt+3\dt^2,\bar t +\dt^2]$},
\end{cases}
\eqnn
where $\fff{t_1,t_2}$ is the flow generated by $v$. Keep $ h^+$. We use the notation $\LL_{\bar t}$ for the solution $\hat \mu$ of \r{pde} with $\hat v$ and $\hat h$, i.e. $\LL_{\bar t}( \mu):=\hat \mu$. We also denote $\LL_{\bar t}(\mu, v, h):=(\hat \mu,\hat v,\hat h)$.\\
\b{Transformation RIGHT $\RR$}: The idea is to replace the increase-transport with the transport-increase. Let $(\mu, v, h)$ be given, and $\bar t$ be a time such that: $ h^-_t=0$ on the interval $[\bar t-\dt^2,\bar t +\dt-2\dt^2]$; $ v_t=0$ on the interval $[\bar t-\dt^2,\bar t]$; $ h^+_t=0$ on the interval $[\bar t,\bar t+\dt-2\dt^2]$. Then replace $ v$ with $\hat v$ defined as follows:
\bqn
\hat v_t:=\begin{cases}
 v_t&\mbox{~~for $t\in[0,\bar t-\dt^2]\cup(\bar t +\dt-2\dt^2,1]$},\\
 v_{t+\dt^2}&\mbox{~~for $t\in(\bar t-\dt^2,\bar t+\dt-3\dt^2]$},\\
0&\mbox{~~for $t\in(\bar t+\dt-3\dt^2,\bar t +\dt-2\dt^2]$}.
\end{cases}
\eqnn
Also replace $ h^+$ with $\hat h^+$ defined as follows:
\bqn
\hat h^+_t:=\begin{cases}
 h^+_t&\mbox{~~for $t\in[0,\bar t-\dt^2]\cup(\bar t +\dt-2\dt^2,1]$},\\
0&\mbox{~~for $t\in(\bar t-\dt^2,\bar t+\dt-3\dt^2]$},\\
\fff{\bar t,\bar t+\dt-2\dt^2}\# h^+_{t-\dt+2 \dt^2}&\mbox{~~for $t\in(\bar t+\dt-3\dt^2,\bar t +\dt-2\dt^2]$}.
\end{cases}
\eqnn
Keep $ h^-$. We use the notation $\RR_{\bar t}$ for the solution $\hat \mu$ of \r{pde} with $\hat v$ and $\hat h$, i.e. $\RR_{\bar t}( \mu):=\hat \mu$. We also denote $\RR_{\bar t}(\mu, v, h):=(\hat \mu,\hat v,\hat h)$.

\newcommand{\RLD}{\RR\LL\DD}
We define $\DD$ as the composition $\DD:=\DD_{\bar t_n}\circ\DD_{\bar t_{n-1}}\circ\ldots\circ \DD_{\bar t_2}\circ\DD_{\bar t_1}$ where $\bar t_1<\bar t_2<\ldots<\bar t_n$ are all times in the set $\Pg{0,\dt^2,2\dt^2,\ldots,(2^{2k}-1)\dt^2,1}$ such that $\DD_{\bar t}$ can be applied. We define $\LL,\RR$ similarly. Finally, we define $\RLD$ as the composition $\RR\circ\LL\circ\DD$. We apply $\RLD$ iteratively to $\bmu$. One can observe that, after $2^{k}-1$ iterations, the result is a fixed point for $\RLD$, i.e. $\RLD(\RLD^{(2^k-1)}(\bmu))=\RLD^{(2^k-1)}(\bmu)$. We define $\bbmu:=\RLD^{(2^k-1)}(\bmu)$ such fixed point.

One can observe that $\bbmu$ is the solution of \r{pde} for a certain $\bbv,\bbh$ (depending on $v^k,h^k$) of this kind:
$$\bbv_t=0 \mbox{~~for $t\in[0,\dt]\cup (1-\dt,1]$},\qquad\qquad\qquad \bbh_t=\begin{cases}
-(\tilde h^k_t)^-&\mbox{~~for $t\in[0,\dt]$},\\
0 &\mbox{~~for $t\in (\dt,1-\dt]$},\\
(\tilde h^k_t)^+&\mbox{~~for $t\in(1-\dt,1]$}.
\end{cases}$$

\b{Step 2.2:} We now prove three properties of $\bbmu$:\\
1. $\bbmu$ drives $\mu_0$ to $\bmu_1$, i.e. $(\bbmu,\bbv,\bbh)\in V(\mu_0,\bmu_1)$. Indeed, transformations $\DD,\LL,\RR$ do not change initial and final times.\\
2. It holds $\fb{\bbmu,\bbv,\bbh}\leq\fb{\bmu,v^k,h^k}$. Indeed, it is easy to prove the following properties
$$\fb{\DD_\bart(\mu, v, h)}= \fb{\mu, v, h},\qquad
\fb{\LL_\bart(\mu, v, h)}\leq \fb{\mu, v, h},\qquad
\fb{\RR_\bart(\mu, v, h)}\leq \fb{\mu, v, h}.$$
3. It holds $\tw{\mu_0,\bmu_1}\leq \fb{\bbmu,\bbv,\bbh}$. Observing the explicit structure of $\bbmu$ in which one has remove of mass in $[0,\dt]$, then transport in $[\dt,1-\dt]$, then creation of mass in $[1-\dt,1]$ one can take $\bbmu_\dt\sotto\bbmu_0=\mu_0$, $\bbmu_{1-\dt}\sotto\bbmu_1=\bmu_1$ and $\bbmu_{1-\dt}=\tilde\Phi^k_{[\dt,1-\dt]}\#\bbmu_\dt$ to estimate
\bqn
&&\tw{\mu_0,\bmu_1}\leq a^2\Pt{|\bbmu_0-\bbmu_\dt|+|\bbmu_1-\bbmu_{1-\dt}|}^2+b^2W_2^2(\bbmu_\dt,\bbmu_{1-\dt}).
\eqnl{e-dai}
Using the standard Benamou-Brenier formula \r{e-bb} for the last term and the change of variable $\tau\to t=(1-2\dt)\tau+\dt$, we have
$$W_2^2\Pt{\bbmu_\dt,\bbmu_{1-\dt}}\leq (1-2\dt) \int_0^{1}  dt \Pt{\int_{\R^d}d \bbmu_t \,|\tilde v^k_t|^2}\leq \int_0^{1}  dt \Pt{\int_{\R^d}d \bbmu_t \,|\tilde v^k_t|^2},$$
that, applied to \r{e-dai}, gives $\tw{\mu_0,\bmu_1}\leq \fb{\bbmu,\bbv,\bbh}$.

\b{Step 2.3:} We now prove $\tw{\mu_0,\mu_1}\leq\fb{\mu,v,h}$. For each $k$ it holds $\tw{\mu_0,\bmu_1}\leq \fb{\bbmu,\bbv,\bbh}\leq \fb{\bmu,v^k,h^k}$. Since $\lim_k |\gw{\mu_0,\bmu_1}-\gw{\mu_0,\mu_1}|\leq \lim_k \gw{\bmu_1,\mu_1}\leq \lim_k d(\bmu,\mu)=0$, then $\lim_k \tw{\mu_0,\bmu_1}=\tw{\mu_0,\mu_1}$. Then 
\bqn
\tw{\mu_0,\mu_1}=\lim_k \tw{\mu_0,\bmu_1}\leq\lim_k \fb{\bmu,v^k,h^k}=\fb{\mu,v,h}.
\eqnn

\b{Step 3.} We now prove Theorem \ref{t-gbb}. We divide the proof in two parts. In part 3.1, we generalize the inequality $\tw{\mu_0,\mu_1}\leq \fb{\mu,v,h}$. In Part 2, we prove the converse inequality.

\b{Step 3.1.} We first prove that $\tw{\mu_0,\mu_1}\leq \fb{\mu,v,h}$, with less regularity requirement. For $h$, we pass from $L^\infty$ to $L^1$ regularity. On the side of $v$, we pass from Lipschitz continuity with respect to space and uniform boundedness to $v_t\in L^2(dt\, d\mu_t)$.

First, one can easily pass from the case of $h$ in $L^\infty$ to the case of $h$ in $L^1$. The idea is to define $\bmu$ as in Step 1, and to provide similar estimates. Instead of a global constant $P$, one needs to define $$p_n^k:=\int_{n2^{-k}}^{(n+1)2^{-k}}dt|h_t|,$$ then prove 
\bqn
\gw{\bmu_{n\dt},\bmu_{n\dt+\tau}}\leq 2ap_n^k+ bMm\dt.
\eqnl{e-quasilip1}
and
\bqn
\gw{\bmu[k-1]_{(n+2) \dt},\bmu_{(n+2) \dt}}&\leq& e^{ L\dt}\gw{\bmu[k-1]_{n\dt},\bmu_{n\dt}}+be^{ L\dt}\dt M  p_n^k +b\dt M p_n^k.
\eqnn
This implies
$$\gw{\bmu[k-1]_{n \dt},\bmu_{n \dt}}\leq 3 b M e^{ L/2} \Pt{\int_0^1 dt |h_t|}  2^{-k},$$
hence, summing up, we have
\bqn
\ww{\bmu[k],\bmu[k+l]}\leq 4 a\psi(2^{-k+1})+C 2^{-k}
\eqnn
with $C_2:= 4bMm+6 b M e^{ L/2} \Pt{\int_0^1 dt |h_t|}$ and $\psi(\eps):=\sup_{t\in[0,1-\eps]} \int_t^{t+\eps} |h_t|$ that satisfies $\psi(\eps)\rightarrow 0$ for $\eps\rightarrow 0$. Hence $\bmu$ is a Cauchy sequence in $C([0,1],\Mu)$.

The proof that the limit $\mu^*=\lim_k \bmu$ coincides with $\mu$ is equivalent to Part 1.3. Finally, one can easily prove $\fb{\mu^k,v^k,h^k}\rightarrow \fb{\mu,v,h}$ by following the estimates of Part 1.4.

We now generalize our result to $v\in L^2(dt\,d\mu_t)$. The proof is completely equivalent to the generalization of the proof of the Benamou-Brenier formula given in \cite[Theorem 8.1]{villani}, Step 2. The main idea is to introduce the variable $m_t:=\rho_t v_t$, where $\rho_t$ is the density of $\mu_t$,  and observe that $\rho_t |v_t|^2=|m_t|^2/\rho_t$ is a convex function of $\rho_t,m_t$. Then, we write $\fb{\mu, m, h}$ with an abuse of notation, and observe that it is convex with respect to its arguments. The presence of the term $h$ makes no difference on this point with respect to \cite[Theorem 8.1]{villani}, Step 2.

Summing up, we have $$\tw{\mu_0,\mu_1}=\lim_{\lam\to 0}\tw{\mu^\lam_0,\mu^\lam_1}\leq \lim_{\lam\to 0} \fb{\mu^\lam,m^\lam,h}\leq \fb{\mu,m,h}$$
with $h\in L^1(dt\, d\mu_t)$ and $v_t\in L^2(dt\, d\mu_t)$.

\b{Step 3.2.}  We now prove that $\inf\Pg{\fb{\mu,v,h}\ |\ (\mu,v,h)\in V(\mu_0,\mu_1)}\leq \tw{\mu_0,\mu_1}$ by giving a sequence $(\mu^k,v^k,h^k)$ realizing the equality at the limit. First of all, observe that there exists\footnote{The result can be proven even without assuming the existence of $\tilde\mu_0,\tilde\mu_1$, via a double limit and a diagonalization argument.} a choice $\tilde\mu_0,\tilde\mu_1$ such that
\bqn
\tw{\mu_0,\mu_1}=a^2\Pt{|\mu_0-\tilde\mu_0|+|\mu_1-\tilde\mu_1|}^2 +b^2 W_2^2(\tilde\mu_0,\tilde\mu_1),
\eqnn
and with $\tilde\mu_0\sotto \mu_0,\tilde\mu_1\sotto\nu_1$. Define  $\psi$ to be the optimal map realizing $W_2^2(\tilde\mu_0,\tilde\mu_1)$, that exists since $\tilde\mu_0,\tilde\mu_1\in\Muu$. Also define (see \cite{villani})
$$\psi_t(x):=(1-t)x+t\psi(x),\qquad v^*_t:=(\psi-\Id)\circ \psi_t^{-1},\qquad \tilde\mu_t:=\psi_t\#\tilde\mu_0,$$
and recall that $(\tilde\mu,v^*)$ is the choice realizing the equality in the standard Benamou-Brenier formula \r{e-bb}, i.e.
$$W_2^2(\tilde\mu_0,\tilde\mu_1)=\fa{\tilde\mu,v^*}.$$

Then, write a dynamics first driving $\mu_0$ to $\tilde\mu_0$ via removal of mass, then $\tilde\mu_0$ to $\tilde\mu_1$ via push-forward of measure, and finally $\tilde\mu_1$ to $\mu_1$ with creation of mass. More precisely, fix an integer $k$, $\dt:=2^{-k}$ and define $v^k,h^k$ as follows:
\bqn
v^k_t:=\begin{cases}
0&\mbox{~~~for $t\in[0,\dt]\cup (1-\dt,1]$},\\
(1-2\dt)^{-1}v^*_{\frac{t-\dt}{1-2\dt}} & \mbox{~~~for $t\in(\dt,1-\dt]$},
\end{cases}
\quad
h^k_t:=\begin{cases}
-\dt^{-1}(\mu_0-\tilde\mu_0)&\mbox{~~~for $t\in[0,\dt]$},\\
0 & \mbox{~~~for $t\in(\dt,1-\dt]$},\\
\dt^{-1}(\mu_1-\tilde\mu_1)&\mbox{~~~for $t\in[1-\dt,1]$}.\\
\end{cases}
\eqnn
The corresponding solution $\mu^k$ of \r{pde} with vector field $v^k$ and source $h^k$ satisfies $(\mu^k,v^k,h^k)\in V(\mu_0,\mu_1)$ and
\bqn
\int_0^1 dt\Pt{\int_{\R^d} d |h_k|}=|\mu_0-\tilde\mu_0|+|\mu_1-\tilde\mu_1|,\ W_2^2(\tilde\mu_0,\tilde\mu_1)=\int_0^1 d \tau \Pt{\int_{\R^d}d\tilde\mu_\tau |v^*_\tau|^2}=(1-2\dt) \int_0^1 dt \Pt{\int_{\R^d}d\mu_t |v^k_\tau|^2}.
\eqnn
One then has 
\bqn
\fb{\mu^k,v^k,h^k}=a^2\Pt{|\mu_0-\tilde\mu_0|+|\mu_1-\tilde\mu_1|}^2+ b^2 (1-2\dt)^{-1} W_2^2(\tilde\mu_0,\tilde\mu_1)\leq (1-2^{-k+1})^{-1} \tw{\mu_0,\mu_1}.
\eqnn
Passing to the limit, we have the result
\bqn
\inf\Pg{\fb{\mu,v,h}\ |\ (\mu,v,h)\in V(\mu_0,\mu_1)}\leq \lim_k \fb{\mu^k,v^k,h^k}\leq \tw{\mu_0,\mu_1}.
\eqnn
\eproof

\noindent \b{Acknowledgments}: The authors thank Luigi Ambrosio for suggesting looking for a generalized Benamou-Brenier formula and
 acknowledge the support of the NSF Grant \#1107444 (KI-Net).\\
This work was partly funded by Carnot STAR Institute in the framework of a researcher exchange program. It was conducted during a visit of F. Rossi to Rutgers University, Camden, NJ, USA. F. Rossi also thanks the institution for its hospitality.

\end{document}